\documentclass{article}
\usepackage{amssymb,amsmath}

\newcommand{\dps}{\displaystyle}
\newcommand{\ptl}{\partial}
\newcommand{\ups}{\upsilon}
\newcommand{\eps}{\epsilon}
\newcommand{\ra}{\rightarrow}
\newcommand{\lra}{\longrightarrow}
\newcommand{\g}{\gamma}
\newcommand{\vp}{\varphi}
\newcommand{\om}{\Omega}
\newcommand{\var}{\varepsilon}
\newcommand{\dt}{\delta}
\newcommand{\te}{\theta}
\newcommand{\na}{\nabla}
\newcommand{\tna}{\tilde{\na}}
\newcommand{\vy}{1+\var-y}
\newcommand{\tC}{\tilde{C}}
\newcommand{\G}{\Gamma}
\newcommand{\brx}{\bar{\mathrm{x}}}
\newcommand{\rx}{\mathrm{x}}

\newtheorem{theorem}{Theorem}

\newtheorem{remark}[theorem]{Definition}

\newtheorem{lemma}[theorem]{Lemma}

\newtheorem{proposition}[theorem]{Proposition}

\newtheorem{corollary}[theorem]{Corollary}

\begin{document}

\title{\bf Global Well-posedness and Regularity of
Weak Solutions to the Prandtl's System}

\author{{Zhouping Xin\thanks{The Institute of Mathematical Sciences \& Department of
Mathematics, The Chinese University of Hong Kong, Shatin, N.T.,
Hong Kong. Email: zpxin@ims.cuhk.edu.hk. This research was supported in part by grants
from RGC of HKSAR 14300917, 14300819, 14302819, 14301421, Basic and Applied Basic Research Foundation of Guangdong Province 2020131515310002, and the Key Project of National Nature Science Foundation of China (Grant No. 12131010).}} \vspace{3mm} \\
{Liqun Zhang\thanks{Institute of Mathematics, Academy of Mathematics and System Science,
Chinese Academy of Sciences, Beijing 100080, China. Email: lqzhang@math.ac.cn.  This research was supported in part by NSFC grants No.10325104, 11471320 and 11631008}} \vspace{3mm} \\
{Junning Zhao\thanks{Department of Mathematics, Xiamen University,
Xiamen, Fujian 361005, China. Email: jnzhao@jingxian.xmu.edu.cn}
}}

\date{}

\maketitle

\begin{center}
\begin{minipage}{13cm}
\begin{center}
{\bf Abstract}
\end{center}
We continue our study on the global solution to the
two-dimensional Prandtl's system for unsteady boundary layers in
the class considered by Oleinik provided that the pressure is
favorable.  First, by using a different method from [13], we gave a
direct proof of existence of a global weak solution by a direct
BV estimate.  Then we prove the uniqueness and continuous dependence on data of such a weak solution to the
initial-boundary value problem.  Finally, we
show the smoothness of the weak solutions and then the global existence of smooth solutions.
\end{minipage}
\end{center}

\vskip 2cm

\section{Introduction}
\indent

\setcounter{equation}{0}

In this paper, we continue our study on the initial-boundary value
problem for the two-dimensional unsteady Prandtl's system:
\begin{eqnarray}
\left\{
\begin{array}{ll}
\ptl_t u + u \ptl_x u + \ups \ptl_y u + \ptl_x P = \nu \ptl^2_x u, \qquad 0<x<L, \qquad y>0 \\
\ptl_x u + \ptl_y \ups = 0 \\
u|_{t=0} = u_0 (x,y), \qquad u|_{y=0} = 0 \\
\ups|_{y=0} = \ups_0 (x,t), \qquad u|_{x=0} = u_1 (y,t) \\
u(x,y,t) \ra U(x,t), \qquad y \ra + \infty
\end{array}
\right.
\end{eqnarray}
with $\nu$ being a fixed positive constant, and the pressure $p$
determined by the Bernoulli's law:
\begin{eqnarray}
\ptl_t U + U \ptl_x U + \ptl_x P = 0.
\end{eqnarray}

The physical situation described by problem (1.1) with (1.2)
corresponds to a plane unsteady flow of viscous incompressible
fluid in the presence of an arbitrary injection and removal of
fluid across the boundaries.  Thus, one may assume that
\begin{eqnarray}
U(x,t) > 0, \quad u_0 (x,t) > 0, \quad u_1 (y,t) > 0, \quad \text{and} \quad v_0 (x,t) \leq 0
\end{eqnarray}

As in Oleinik [8, 9], we assume that the data are in the monotone
class in the sense that
\begin{eqnarray}
\ptl_y u_0 (x,y) > 0, \quad \ptl_y u_1 (t,y) > 0
\end{eqnarray}

Under the further assumption that the pressure is favorable, i.e.,
\begin{eqnarray}
\ptl_x P (x,t) \leq 0 \quad {\rm for} \quad t>0, \quad 0 < x \leq
L,
\end{eqnarray}
we have shown the existence of a global weak solution to the
Cauchy problem (1.1) by a splitting method in [13].  The main
purpose in this paper is to show that such a weak solution is in
fact unique and depend continuously on the data.  Furthermore,
such a solution in fact is a classical solution in the sense that it is smooth in the interior.

To formulate the problem, as in [8], we introduce the following
Crocco transformation:
\begin{eqnarray}
\tau = t, \quad \xi = x, \quad \eta = \frac{u(x,y,t)}{U(x,t)},
\quad w(\tau, \xi, \eta) = \frac{\ptl_y u(x,y,t)}{U(x,t)} \, .
\end{eqnarray}

Then the initial-boundary value problem (1.1) is transformed into
the following initial-boundary value problem:
\begin{eqnarray}
\left\{
\begin{array}{ll}
\ptl_\tau w^{-1} + \eta U \ptl_\xi w^{-1} + A \ptl_\eta w^{-1} -
B w^{-1} = - \nu \ptl^2_\eta w \\
\qquad \qquad \qquad \qquad {\rm on} \quad Q=\{ ( \xi, \eta,
\tau), \quad 0<\tau<\infty, \quad 0<\xi<L, \quad 0<\eta<1 \}
\\
\dps w|_{\tau=0} = \frac{\ptl_y u_0}{U} \equiv w_0, \quad
w|_{\eta=1} =
0 \\
\dps w|_{\xi=0} = w_1 \equiv \frac{\ptl_y u_1 (y,\tau)}{U(0,\tau)},
\quad (\nu w \ptl_\eta w - v_0 w)|_{\eta=0} = \frac{\ptl_x P}{U},
\end{array}
\right.
\end{eqnarray}
where
\begin{eqnarray}
\left\{
\begin{array}{ll}
\dps A= (1-\eta^2) \ptl_x U + (1-\eta) \frac{\ptl_t U}{U}, & {\rm and}\\
\dps B= \eta \ptl_x U + \frac{\ptl_t U}{U} .
\end{array}
\right.
\end{eqnarray}

Set
\begin{eqnarray*}
Q_T = \{ (\xi, \eta, \tau) | \quad 0<\tau<T, \quad 0<\xi<L, \quad
0<\eta<1 \}
\end{eqnarray*}

\vskip 1cm

A weak solution to the initial-boundary value
problem (1.7) can be defined as follows:

\vskip 5mm

\begin{remark}
{\rm A function $w \in BV (Q_T) \cap L^\infty (Q_T)$ is said to
be a weak solution to the problem (1.7) if the following
conditions are satisfied:
\begin{enumerate}
    \item[(i)] There exists a positive constant $C=C(T)$ such that
    \begin{eqnarray}
    C^{-1} (1-\eta) \leq w(\tau, \xi, \eta) \leq C(1-\eta), \quad
    \forall(\xi, \eta, \tau) \in Q_T, \quad {\rm and}
    \end{eqnarray}
    \begin{eqnarray}
    (1-\eta)^{\frac{1}{2}} \, \ptl_\eta w \in L^2 (Q_T) ,
    \end{eqnarray}

    \item[(ii)] $w_{\eta \eta}$ is a locally bounded measure in
    $Q_T$, and
    \begin{eqnarray}
    \int\!\!\!\int_{Q_T} (1-\eta)^2 d|w_{\eta \eta}| < \infty ,
    \end{eqnarray}
    where $d|w_{\eta \eta}|$ denotes the variation of $w_{\eta
    \eta}$.

    \item[(iii)] The boundary conditions in (1.7) are satisfied in
    the sense of trace, i.e.,
    \begin{eqnarray}
    \left\{
    \begin{array}{lll}
    \g w (\xi, \eta, \tau)|_{\tau=0} = w_0 (\xi, \eta) & {\rm
    a.e.\ on} & \bar{Q}_T \cap \{ \tau=0\} , \\
    \g w (\xi, \eta, \tau)|_{\eta=1} = 0 & {\rm
    a.e.\ on} & \bar{Q}_T \cap \{ \eta=1\} , \\
    \g w (\xi, \eta, \tau)|_{\xi=0} = w_1 (\eta, \tau) & {\rm
    a.e.\ on} & \bar{Q}_T \cap \{ \xi=0\} , \\
    \dps \g w_\eta (\xi, \eta, \tau)|_{\eta=0} = w_0 (\xi, \tau) +
    \left. \frac{\ptl_x P}{U} \frac{1}{\g w} \right|_{\eta = 0} & {\rm
    a.e.\ on} & \bar{Q}_T \cap \{ \eta=0\} , \\
    \end{array}
    \right.
    \end{eqnarray}
    where $\g w$ and $\g w_\eta$ are the traces of $w$ and
    $w_\eta$ on the corresponding boundaries respectively.

    \item[(iv)] For any $\vp \in C^1 (\bar{Q}_T)$ with $\vp|_{\tau=0} = \vp|_{\xi=0} = \vp|_{\xi=L} =
    0$, the following identity holds:
    \begin{eqnarray}
    \begin{array}{ll}
    & \dps - \int_\om w^{-1} \vp (1-\eta)^2 \, d\xi \, d\eta|_{\tau=\tau} +
    \int\!\!\!\int_{Q_T} w^{-1} (1-\eta)^2 \, \ptl_\tau \vp \, d\xi \,
    d\eta \, d\tau \\
    + & \dps \nu \int\!\!\!\int_{Q_T} \left( (1-\eta)^2 \vp \right)_\eta w_\eta \, d\xi \,
    d\eta \, d\tau + \int\!\!\!\int_{Q_T} (\eta U \vp)_\xi (1-\eta)^2 \, w^{-1} \, d\xi \,
    d\eta \, d\tau \\
    + & \dps \int\!\!\!\int_{Q_T} \left( (1-\eta)^2 A \vp \right)_\eta w^{-1} \, d\xi \,
    d\eta \, d\tau + \int\!\!\!\int_{Q_T} (1-\eta)^2 \vp B w^{-1} \, d\xi \,
    d\eta \, d\tau \\
    + & \dps \int^\tau_0\!\!\!\!\int^L_0 \ups_0 (\xi, \tau) \vp (\xi, 0, \tau) d\xi \,
    d\tau = 0 \qquad \qquad \forall \, \tau \in (0,T) \, .
    \end{array}
    \end{eqnarray}
\end{enumerate}
}
\end{remark}

\vskip 1cm

The main results in this paper can be summarized in the following
theorem:

\vskip 5mm

\begin{theorem}
{\rm Assume that the data satisfy condition (1.3) and (1.4), and
that the pressure is favorable, i.e., (1.5) holds.  Then,
\begin{enumerate}
    \item[(i)] There exists a weak solution $w \in BV (Q_T) \cap L^\infty
    (Q_T)$ in the sense of Definition 1.1;
    \item[(ii)] Such a weak solution is unique and depends
    continuously on the initial and boundary data in $L^1$-norms;
    \item[(iii)] Such a weak solution is smooth in $Q_T$ for any $T>0$.
\end{enumerate}
 }
\end{theorem}

{\bf Remark:} It should be remarked that the existence of a weak solution has
been proved by the first two authors in [13], in this paper, we
give a different proof by a direct BV estimate, which yields also some additional estimates (1.10) and (1.11) that are important for the structure of such weak solutions. The key new ingredients of the current paper are the uniqueness and the regularity of a weak solution.
Some new ideas are required due to the strong degeneracy of the
equation in (1.7). Moreover, the well-posedness we proved in the BV class in Proposition 3.5, in particular the estimate (3.35) is of independent interests and optimal, and the Poincar\'e type inequality we introduced for the proof of H\"older regularity of weak solutions is also of independent interests.

\vskip 5mm
Together with Proposition 3.5, we have
\begin{corollary}
{\rm Assume that there are two initial datum $u_{10}(x,y)$, $u_{20}(x,y)$, boundary datum  $u_{11}(y,t)$, $u_{21}(y,t)$ and $v_{10}(x,t)$, $v_{20}(x,t)$ satisfying condition (1.3) and (1.4) respectively, with the same $U(x,t)$ as $y\rightarrow \infty$ and
that the pressure is favorable, namely, (1.5) holds. Let $u_1(t,x,y)$ and $u_2(t,x,y)$ be the corresponding solution to the problem (1.1). Then, there exists a constant $C=C(T,L)$ such that
\begin{eqnarray}
\begin{array}{lcl}
& & \dps \int^L_0\!\!\!\int^{\infty}_0 |\partial_y u_1(x,y,t) -\partial_y u_2(x,\tilde y,t)| \frac{\partial_y u_1(x,y,t)}{U^2(x,t)} dy\,dx \\
& \leq & \dps C \left\{ \int^{\infty}_0\!\!\!\int^L_0 |\partial_y u_{10} (x,y) -\partial_y
u_{20}(x,\tilde y)|\frac{\partial_y u_{10}(x,y)}{U^2(x, 0)} dx\,dy + \int^t_0\!\!\!\int^L_0 |\ups_{10}
(x,s) - \ups_{20}(x,s)| dx\,ds
\right. \\
& &  \dps + \left. \int^t_0\!\!\!\int^{\infty}_0 |\partial_y u_{11} (y,s)
-\partial_y u_{21}(\tilde y,s)|\frac{\partial_y u_{11}(y,s)}{U^2(x,s)} dy\,ds   \right\} \,,
\end{array}
\end{eqnarray}
}
where $\tilde y$ is given by $u_1(x,y,t)=u_2(x,\tilde y,t)$ which is uniquely determined.
\end{corollary}

\vskip 1cm

Finally, we mention that there have been increasing activities in recent years on the studies of the unsteady Prandtl's system in both two and three dimensions with substantial results for general initial data, ill-posedness, instability and different methods, we refer to [22-32] and references therein for these further developments.

\newpage

\section{Existence}
\indent

\setcounter{equation}{0}

In this section, we will give an alternative proof of existence of
weak solutions to the initial-boundary value problem (1.7) by a
vanishing viscosity method, which is somewhat simpler than the
approach given by the first two authors in [13].

For the convenience of presentation, we may rewrite the problem
(1.7) as (we use notations independent of those in \S 1).

\begin{eqnarray}
\left\{
\begin{array}{ll}
\ptl_t u - u^2 \, \ptl^2_y \, u + a \, \ptl_x \, u + b \, \ptl_y \, u + cu = 0 \\
u(x,y,t)|_{t=0} = u_0 (x,y) \\
u(x,y,t)|_{x=0} = u_1 (y,t) \\
\dps u(x,y,t)|_{y=1} = 0, \quad u \, \ptl_y \, u|_{y=0} = \left. \left(
\ups_0 \, u + \frac{\ptl_x P}{U} \right) \right|_{y=0}
\end{array}
\right.
\end{eqnarray}
with
\begin{eqnarray}
\left\{
\begin{array}{ll}
\dps a(x,y,t)=y U(x,t), \quad b(x,y,t) = A = (1-y^2) \ptl_x U(x,t)
+ (1-y) \frac{\ptl_t U}{U} \, , \\
\dps c=c(x,y,t) = B = (1-y) \ptl_x U(x,t) - \frac{\ptl_x P}{U} \,
.
\end{array}
\right.
\end{eqnarray}

We also use the notations
\begin{eqnarray*}
\begin{array}{ccl}
Q_T & = & \{ (x,y,t) | 0<x<L, \quad 0<y<1, \quad 0<t<T \} , \\
\om & = & \{ (x,y) | 0<x<L, \quad 0<y<1 \} .
\end{array}
\end{eqnarray*}

\vskip 1cm

\setcounter{theorem}{0}

Then our main results in this section are:

\vskip 5mm

\begin{theorem}
{\rm Assume that the pressure is favorable, i.e., (1.9) holds, and
there exists a positive constant $C_0$ such that for all
$0<x<L $, \, $0<y<1$, \, $0<t<T$, the initial and boundary condition
\begin{eqnarray}
C^{-1}_0 (1-y) < u_0 (x,y) < C_0 (1-y), \quad C^{-1}_0 (1-y) < u_1
(y,t) < C_0 (1-y) .
\end{eqnarray}

Then there exists a weak solution to the initial-boundary value
problem (2.1) in the sense given in Definition 1.1, i.e., there
exists a $u \in BV(Q_T) \cap L^\infty (Q_T)$ with the properties:
\begin{enumerate}
    \item[(i)] It holds that
    \begin{eqnarray}
    (1-y)^{\frac{\alpha - 1}{2}} \, \ptl_y \, u \in L^2 (Q_T), \quad
    {\rm and}
    \end{eqnarray}
    \begin{eqnarray}
    C^{-1}_1 (1-y) \leq u(x,y,t) \leq C_1(1-y), \quad (x,y,t) \in
    \, Q_T
    \end{eqnarray}
    for some positive constants $C_1 = C_1(T,\om)$ and $\alpha
    >0$.

    \item[(ii)] $u_{yy}$ is a locally bounded measure in $Q_T$ in
    the sense that
    \begin{eqnarray}
    \int\!\!\!\int_{Q_T} (1-y)^\alpha \, d|u_{yy}| < + \infty ,
    \end{eqnarray}
    where $d|u_{yy}|$ denotes the variation of $u_{yy}$, and $\alpha >
    0$ is a constant.

    \item[(iii)] And
    \begin{eqnarray}
    \left\{
    \begin{array}{lll}
    \g u (x,y,t)|_{t=0} = u_0 (x,y) & {\rm
    a.e.\ on} & \om , \\
    \g u (x,y,t)|_{y=1} = 0 & {\rm
    a.e.\ on} & \bar{Q}_T \cap \{ y=1 \} , \\
    \g u (x,y,t)|_{x=0} = u_1 (y,t) & {\rm
    a.e.\ on} & \bar{Q}_T \cap \{ x=0 \} , \\
    \dps \g u_y (x,y,t)|_{y=0} = \left. \left( \ups_0 +
    \frac{\ptl_x P}{U} \frac{1}{\g u} \right) \right|_{y = 0} & {\rm
    a.e.\ on} & \bar{Q}_T \cap \{ y=0 \} , \\
    \end{array}
    \right.
    \end{eqnarray}
    where $\g u$ and $\g u_y$ are the traces of $u$ and $\ptl_y \,
    u$ on the corresponding boundaries respectively.

    \item[(iv)] For any $\vp \in C^1 (\bar{Q}_T)$ with $\vp|_{t=0} =
    \vp|_{x=0} = \vp|_{x=L} = 0$, the following identity holds:
    \begin{eqnarray}
    \begin{array}{ll}
    & \dps - \int_\om u^{-1} \vp (1-y)^2 dx\,dy|_{t=t} +
    \int\!\!\!\int_{Q_T}[u^{-1} (1-y)^2 \, \ptl_t \vp  +
     \left( (1-y)^2 \vp \right)_y \,
    u_y]\,dx\,dy\,dt \\
    + & \dps \int\!\!\!\int_{Q_T} [(a \vp)_x (1-y)^2 \frac{1}{u} \,
     + \left( (1-y)^2 b \vp \right)_y \, \frac{1}{u} \,
     + (1-y)^2 \, \vp c \, \frac{1}{u} ]\,
    dx\,dy\,dt \\
    + & \dps \int^t_0\!\!\int^L_0 \ups_0 (x,t) \, \vp (x,y,t)
    dx\,dt|_{y=0} = 0 \qquad \qquad \qquad {\rm for\ any} \qquad t \in (0,T) \,
    .
    \end{array}
    \end{eqnarray}
\end{enumerate}
}
\end{theorem}

We will prove the above theorem by studying the following
regularized problem:
\begin{eqnarray}
\left\{
\begin{array}{ll}
\ptl_t \, u - (u+\var)^2 \, \ptl^2_y \, u + (a+\var) \ptl_x \, u +
b \, \ptl_y \, u + c u = 0 \\
u|_{t=0} = u^\var_0, \quad u|_{x=0} = u^\var_1 \\
\dps u^\var |_{y=1} = 0, \quad (u + \var) \ptl_y u|_{y=0} = \left.
\left( \ups^\var_0 (u + \var) + \frac{\ptl_x P}{U} \right)
\right|_{y=0}
\end{array}
\right.
\end{eqnarray}
where $u^\var_0$, $u^\var_1$, and $\ups^\var_0$ are some suitable
regularization of $u_0$, $u_1$, and $\ups_0$ respectively, and $\varepsilon\in(0,1)$ is a constant.

It follows from the standard theory for initial-boundary value
problems for ultra-parabolic equations that for each fixed $\var >
0$, there exists a unique classical solution $u^\var$ to the
problem (2.9).  Our main strategy is to pass the limit as $\var
\ra 0^+$ in (2.9) to obtain a solution to the problem (2.1).  To
this end, some a priori estimates are needed.  We start with the
super-norm estimates.

\vskip 1cm

\begin{lemma}
{\rm Assume that (1.5) holds and that there exists a positive
constant $C_0$ such that the initial and boundary data
\begin{eqnarray}
C^{-1}_0 (1-y) \leq u^\var_i \leq C_0 (1-y), \quad {\rm on} \quad
Q_T, \quad i=0,1 \, .
\end{eqnarray}

Then there exists a positive constant $C_1=C_1(Q_T)$ such that the solution
\begin{eqnarray}
C^{-1}_1 (1-y) \leq u^\var (x,y,t) \leq C_1 (1-y), \quad \forall
(x,y,t) \, \in \, Q_T \, .
\end{eqnarray}
}
\end{lemma}

\vskip 1cm

\begin{enumerate}
    \item[]
    \begin{enumerate}
    \item[{\bf Proof:}] First, it follows from the assumptions (1.3), (1.5), (2.10), and
the standard maximum principle argument that
\begin{eqnarray}
u^\var (x,y,t) \geq 0 \qquad \qquad {\rm on} \qquad Q_T \, .
\end{eqnarray}

Next, we estimate the upper bound on $u^\var$.  For any given $\dt
> 0$, consider
$$
\te (x,y,t) = \frac{u^\var}{1-y+\dt}
$$

Then $\te (x,y,t)$ solves the following initial-boundary value
problem:
\begin{eqnarray}
\left\{
\begin{array}{ll}
\dps \ptl_t \te - (u^\var + \var) \te_{yy} + (a+\var) \ptl_x \te +
\left( b-\frac{2(u^\var + \var)^2}{1-y+\dt} \right) \te_y + \left(
c-\frac{b}{1-y+\dt} \right) \te = 0 \\
\dps \te|_{t=0} = \frac{u^\var_0}{1-y+\dt} \, , \quad \te|_{x=0} =
\frac{u^\var_1}{1-y+\dt} \, , \quad \te|_{y=1} = 0 \\
\dps (u^\var + \var) \left. \left[ \ptl_y \te -
\frac{\te}{1-y+\dt} \right] \right|_{y=0} = \left. \left[ \ups_0
\, \frac{u^\var + \var}{1-y+\dt} + \frac{\ptl_x P}{(1-y+\dt) U}
\right] \right|_{y=0}
\end{array}
\right.
\end{eqnarray}

Due to the specific structure of $b$ and $c$ ((2.2)), the
coefficient  $\dps c(x,y,t)-\frac{b(x,y,t)}{1-y+\dt}$ admits a
uniform bounded,
$$
\left| c(x,y,t) - \frac{b(x,y,t)}{1-y+\dt} \right| \leq \tilde{c}_2
$$
where $\tilde{c}_2$ is independent of $\dt$ and $\var$.  Thus,
standard maximum principle argument leads to
\begin{eqnarray}
0 \leq \te(x,y,t) \leq e^{\tilde{c}_2 T} \, \max \{ \te|_{t=0} ,
\te|_{x=0} , \te|_{y=0} , \te|_{y=1} \} \, .
\end{eqnarray}

Since $\te|_{t=0} \leq C_0$, $\te|_{x=0} \leq C_0$, and
$\te|_{y=1} = 0$, it suffices to estimate the $\dps \max_{y=0} \,
\te(x,y,t)$. At the maximum, $\ptl_y \te \leq 0$, so the boundary
condition at $y=0$ yields
$$
\frac{u^\var + \var}{1+\dt} \te \leq |\ups_0| \frac{u^\var +
\var}{1+\dt} + \frac{|\ptl_x P|}{(1+\dt)|U|} \, ,
$$
which implies
$$
\te^2 \leq |\ups_0| (1+\te) + \frac{|\ptl_x P|}{|U|} \, .
$$

Consequently,
$$
\dps \max_{y=0} \, \te(x,y,t) \leq \tilde{c}_3
$$
where $\tilde{c}_3$ is a positive constant independent  of $\var$
and $\dt$.  Thus (2.14) shows that
$$
0 \leq \te(x,y,t) \leq e^{\tilde{c}_2 T} \, \max (c_0,
\tilde{c}_3) = \tilde{c}_4 \, ,
$$
and so
$$
0 \leq u^\var (x,y,t) \leq \tilde{c}_4 (1-y+\dt) \, .
$$

Since $\dt$ is arbitrary, the desired upper bound in
(2.11) follows. Finally, we obtain the lower bound estimate.  Set
\begin{eqnarray}
\left\{
\begin{array}{rcl}
\te(x,y,t) & = & e^{\alpha t} \, u^\var - \mu_0 \, \vp (y)
e^{-\beta t}
, \quad {\rm with} \\
\vp(y) & = & e^{- \g(1-y)^2} (1-y)
\end{array}
\right.
\end{eqnarray}
where $\alpha$, $\beta$, $\g$ are large and $\mu_0$ is small to be
determined.  Define also operator
$$
\mathcal{L}=\ptl_t - (u^\var + \var)^2 \ptl^2_y + (a+\var) \ptl_x
+ b \, \ptl_y + c \, .
$$

Then direct calculations yield
\begin{eqnarray*}
\begin{array}{rcl}
\mathcal{L} \te & = & \alpha \te + \mu_0 \, \varphi (y) e^{-\beta t}
\left[\alpha + \beta - (u^\var+\var) (6\g + 4\g (1-y)^2) -
b\frac{(2\g (y-1)^2 -1)}{(1-y)}
+ c \right] \, , \\
\te|_{t=0} & = & u^\var_0 - \mu_0 \, \varphi (y) , \quad \te|_{y=1} =
0 \, , \\
\te|_{x=0} & = & u^\var_1 \, e^{\alpha t} - \mu_0 \, \varphi (y)
e^{-\beta t} = e^{\alpha t} \left[ u^\var_1 - \mu_0 \,
e^{-(\alpha+\beta)t} \, \varphi (y) \right] \, , \\
\dps \ptl_y \te|_{y=0} & = & \dps e^{\alpha t} \left. \left(
\ups_0 + \frac{\ptl_x P}{(u^\var +\var) U} \right) \right|_{y=0} -
\mu_0 \, e^{-\g} (2\g -1) e^{-\beta t} \, .
\end{array}
\end{eqnarray*}

We now fix $\g$ so that $\g > \frac{1}{2}$.  Thus, $\ptl_y \,
\te|_{y=0} < 0$.
Then one can choose $\mu_0$ so small such that $\te|_{t=0} \geq 0$, and
$$\te|_{x=0} \geq e^{\alpha t} [u^\var_1 - \mu_0 \, \phi (y) ] \geq
0$$ by the assumption (2.10).

Due to the upper bound estimate on $u^\var$, one can choose
$\alpha$ and $\beta$ big enough so that
$$
\alpha + \beta - (u^\var + \var) (6\g+4\g(1-y)^2) - b\frac
{(2\g(y-1)^2 -1)}{(1-y)} + c > 0 \, .
$$

Hence
$$
\mathcal{L} \te - \alpha \te > 0 \qquad \qquad  {\rm in} \qquad
Q_T \, .
$$

It follows by the maximum principle that $\te \geq 0$ on $Q_T$. As
a consequence, we have shown that
$$
u^\var (x,y,t) \geq C^{-1}_1 (1-y)
$$
where $C_1=C_1(Q_T)$ is a positive constant independent of $\var$.
The proof of Lemma 2.2 is completed.
\begin{flushright}
$\square$
\end{flushright}
    \end{enumerate}
\end{enumerate}

Next, we derive the uniform total variation estimate on the
approximation solutions.  For any given scalar function $\psi$, we set
\begin{eqnarray}
\na \psi = {\rm grad\ } \psi = (\ptl_x \psi, \ptl_y \psi, \ptl_t
\psi)^t \, , \qquad {\rm and} \\
\tna \psi = (\ptl_x \psi, \ptl_t \psi)^t \, . \qquad \qquad \qquad
\qquad \quad \quad
\end{eqnarray}

Then, our next key uniform estimate is the following total
variation estimate on $u^\var$.

\vskip 1cm

\begin{proposition}
{\rm There exists a positive constant $C_2=C_2(Q_T)$ independent
of $\var$, such that
\begin{eqnarray}
\dps \int_{Q_T} |\na u^\var (x,y,t)| dx\,dy\,dt \leq C_2 \, .
\end{eqnarray}

This proposition will follow from the following two lemmas
directly. }
\end{proposition}

\vskip 1cm

\begin{lemma}
{\rm For any given constant $\alpha > -1$, there exists a positive
constant $C_3=C_3(\alpha)$ such that
\begin{eqnarray}
\dps \int_{Q_T} (1-y)^\alpha |\ptl_y u^\var| dx\,dy\,dt \leq C_3
\, , \quad {\rm and} \\
\dps \int_{Q_T} (1-y)^\alpha |\ptl_y u^\var|^2 dx\,dy\,dt \leq C_3
\, . \quad \quad \,
\end{eqnarray}
}
\end{lemma}

\vskip 1cm

\begin{enumerate}
    \item[]
    \begin{enumerate}
        \item[{\bf Proof:}] Note that (2.19) follows immediately
        from (2.20) by the H\"{o}lder inequality.  So it suffices to
        show that (2.20) holds.  Setting $g(u)=ln(u+\var) +
        \var(u+\var)^{-1}$, multiplying the equation in (2.9) by
        $(1-y)^\alpha g'(u^\varepsilon)$, and integrating the resulting equation over
        $Q_T$, one can obtain after integration by parts that
        \begin{eqnarray}
        \begin{array}{ll}
        & \dps \left. \int_\om g(u^\var) (1-y)^\alpha \right|^T_0 dx\,dy \\
        = & \dps - \int_{Q_T} (1-y)^\alpha (\ptl_y u^\var)^2 dx\,dy\,dt +
        \alpha \int_{Q_T} (1-y)^{\alpha-1} \, u^\var \, \ptl_y u^\var
        dx\,dy\,dt \\
        & \dps + \int_{Q_T} g(u^\var) \ptl_y \, ((1-y)^\alpha \, b) \, dx\,dy\,dt +
        \int_{Q_T} (1-y)^\alpha \, \ptl_x a \, g(u^\var) dx\,dy\,dt \\
        & \dps - \int_{Q_T} c(1-y)^\alpha \frac{(u^\var)^2}{(u^\var +
        \var)^2} dx\,dy\,dt + \left. \int^T_0\!\!\!\int^L_0 (1-y)^\alpha \,
        u^\var \ptl_y u^\var dx\,dt \right|^1_0 \\
        & \dps - \left. \int^T_0\!\!\!\int^L_0 (1-y)^\alpha \, b \, g(u^\var)
        dx\,dt \right|^1_0 - \left. \int^T_0\!\!\!\int^1_0 (1-y)^\alpha (a+\var) \, g(u^\var)
        dy\,dt \right|^L_0
        \end{array}
        \end{eqnarray}

        First, it holds that
        \begin{eqnarray}
        \begin{array}{ll}
        & \dps \left| \alpha \int_{Q_T} (1-y)^{\alpha-1} \, u^\var \, \ptl_y u^\var
        dx\,dy\,dt \right| \\
        \leq & \dps \frac{1}{2} \int_{Q_T} (1-y)^\alpha (\ptl_y
        u^\var)^2 dx\,dy\,dt + \frac{\alpha^2}{2} \int_{Q_T}
        (1-y)^{\alpha-2} \, (u^\var)^2 dx\,dy\,dt
        \end{array}
        \end{eqnarray}

        Next, due to Lemma 2.2, one has that
        \begin{eqnarray}
        |(1-y)^\alpha g(u^\var)| \leq C \left\{
        (1-y)^\alpha \, (|lu(1-y)| + 1) \right\} \quad {\rm for}
        \quad (x,y,t) \, \in Q_T
        \end{eqnarray}
        with a positive constant $C=C(Q_T)$.  As a consequence,
        one can derive that
        \begin{eqnarray}
        \begin{array}{ll}
        & \dps \left| \int_\om g(u^\var) (1-y)^\alpha dx\,dy \big|^T_0 \, \right| + \left|
        \int_{Q_T} g(u^\var) \, \ptl_y \, (b(1-y)^\alpha) dx\,dy\,dt
        \right| \\
        + & \, \dps \left| \int_{Q_T} (1-y)^\alpha \, \ptl_x a \, g(u^\var)
        dx\,dy\,dt \right| + \big|
        \int_{Q_T} c(1-y)^\alpha \frac{(u^\var)^2}{(u^\var +
        \var)^2} \, dx\,dy\,dt \\
        + & \dps \left| \int^T_0\!\!\!\int^1_0 (1-y)^\alpha
        (a+\var) \, g(u^\var) dy\,dt \big|^L_0 \, \right| \leq C
        \end{array}
        \end{eqnarray}
        with $C=C(Q_T)$ independent of $\var$, where we have used
        the fact that $|b| \leq C(1-y)$.  Finally, it follows from
        the boundary conditions in (2.9) and the structures of $b$
        that
        \begin{eqnarray}
        \begin{array}{ll}
        & \dps \left| \int^T_0\!\!\!\int^L_0 (1-y)^\alpha \, u^\var \, \ptl_y u^\var dx\,dt \big|^1_0 -
        \int^T_0\!\!\!\int^L_0 (1-y)^\alpha \, b \, g(u^\var) \,
        dx\,dt \big|^1_0 \, \right| \\
        = & \dps \,\, \left| \int^T_0\!\!\!\int^L_0 \left( \frac{u^\var}{\var + u^\var} (u^\var + \var)
        \ptl_y \, u^\var \right)
        \big|_0 dx\,dt - \int^T_0\!\!\!\int^L_0 (b \, g(u^\var))
        \big|_0 dx\,dt \right| \\
        \leq & \dps \int^T_0\!\!\!\int^L_0 \left( \frac{u^\var}{u^\var + \var} \, |\ups^\var_0 (u^\var + \var)|
        + \frac{|\ptl_x P|}{U} \right) (x,0,t) dx\,dt + \tilde{C}
        \\
        \leq & C
        \end{array}
        \end{eqnarray}
        with positive constants $C$ and $\bar{C}$ independent of
        $\var$.  Now, the desired estimate (2.20) follows from
        (2.21) - (2.22), and (2.24) - (2.25).  This completes the
        proof of Lemma 2.4.
\begin{flushright}
$\square$
\end{flushright}
    \end{enumerate}
\end{enumerate}

Next, we estimate $\tna u^\var$. In fact, we have

\vskip 1cm

\begin{lemma}
{\rm There exists a positive constant $C_4=C_4(\om_T)$ independent
of $\var$ such that
\begin{eqnarray}
\sup_{0\leq t\leq T} \int^L_0\!\!\!\int^1_0 \big| \tna u^\var
(x,y,t) \big| dx\,dy \leq C_4 \, .
\end{eqnarray}
}
\end{lemma}

\vskip 1cm

\begin{enumerate}
    \item[]
    \begin{enumerate}
        \item[{\bf Proof:}] For any $\dt>0$, set
        \begin{eqnarray}
        S_\dt (\te)= \left\{
        \begin{array}{cc}
        1 , & \te > \dt , \\
        \dps \frac{\te}{\dt} , & |\te| \leq \dt , \\
        -1 , & \te < -\dt ,
        \end{array}
        \right.
        \end{eqnarray}
        and for any $\xi \in \mathbb{R}^2$,
        \begin{eqnarray}
        I_\dt (\xi) = \int^{|\xi|}_0 S_\dt (s) ds \, .
        \end{eqnarray}

        Then it can be checked easily that
        \begin{eqnarray}
        \na^2 I_\dt (\xi) \geq 0 , \quad \ptl_{\xi_i} I_\dt (\xi)
        = S_\dt (|\xi|) \frac{\xi_i}{|\xi|} \, ,
        \end{eqnarray}
        and
        \begin{eqnarray}
        \lim_{\dt \ra 0^+} \xi \cdot \na^2 I_\dt (\xi) = 0 , \quad {\rm except} \quad |\xi|=0 \,
        .
        \end{eqnarray}

        To show (2.26), it is more convenient to rewrite the
        equation in (2.9) in terms of new dependent variable
        \begin{eqnarray}
        V=\frac{1}{u^\var+\var}
        \end{eqnarray}
        which satisfies
        $$
        \ptl_t V - \ptl_y \left( \frac{V_y}{V^2} \right) +
        (a+\var) \ptl_x V + b \, \ptl_y V - \frac{C
        u^\var}{u^\var+\var} V = 0 \, .
        $$

        Thus,
        \begin{eqnarray}
        \ptl_t \tna V - \ptl^2_y \left( \frac{\tna V}{V^2} \right)
        + \tna ((a+\var) V_x) + \tna (b \,V_y) - \tna \left( \frac{C
        u^\var}{u^\var+\var} V \right) = 0 \, .
        \end{eqnarray}

        Taking inner product of (2.32) with $\dps (\vy)^2 S_\dt
        (\tna V) \frac{\tna V}{|\tna V|}$ and
        integrating the resulting equation over $\om$, one can get
        after integration by parts that

        \begin{eqnarray}
        \begin{array}{ll}
        & \dps \frac{d}{dt} \int_\om I_\dt (\tna V) (\vy)^2
        dx\,dy \\
        = & \dps \int_\om [\ptl^2_y\left(\frac{\tna V}{V^2}\right)\cdot\ptl_{\tna V}
        I_\dt (\tna V) (\vy)^2  -  \tna ((a+\var) V_x)
        \cdot \ptl_{\tna V} I_\dt (\tna V) (\vy)^2] dx\,dy \\
        & \dps - \int_\om [\tna (b V_y) \cdot \ptl_{\tna V}
        I_\dt (\tna V) (\vy)^2  +  \tna \left( \frac{C u^\var}{u^\var
        +\var} V \right) \cdot \ptl_{\tna V} I_\dt (\tna V) (\vy)^2]
        dx\,dy\\
        = & \dps \left\{ \int^L_0 \left[ \ptl_y \left( \frac{\tilde{\na} V}{V^2} \right) \ptl_{\tna
        V}\, I_\dt (\tna V) (\vy)^2 - b(\vy)^2 \, \tna V \, \ptl_{\tna V} I_\dt (\tna V)  \right.
        \right. \\
        & \dps \qquad \qquad + \left. \left. \left. \frac{\tna V}{V^2} \ptl_{\tna V} \, I_\dt (\tna V) \, 2
        (\vy) \right] \right|^1_0 dx \right\} - \left\{ \left.
        \int^1_0 (a+\var) I_\dt (\tna V) (\vy)^2
        dy \right|^L_0 \right\} \\
        & \dps + \left\{ - \int_\om \frac{\tna V_y}{V^2}
        \na^2_{\tna V} \, I_\dt (\tna V) (\tna
        V_y)(\vy)^2 dx\,dy \right\} \\
        & \dps + \left\{ \int_\om \left[
        2\left(\frac{\tna V}{V^2}\right) \ptl_{\tna V} \, I_\dt (\tna
        V) - V_x (\tna a) \ptl_{\tna V} \, I_\dt (\tna V)(\vy)^2 + \ptl_x a I_\dt
        (\tna V)(\vy)^2\right. \right. \\
        & \dps \qquad \qquad + \ptl_y (b(\vy)^2) \tna V \ptl_{\tna V} I_\dt
        (\tna V) + \frac{C u^\var}{u^\var + \var} \tna V \cdot
        \ptl_{\tna V} I_\dt (\tna V) (\vy)^2 \\
        & \dps \qquad \qquad + \left. \left. \var(-C) V
        \tna V \cdot \ptl_{\tna V} \, I_\dt (\tna V) (\vy)^2
        \right] dx\,dy \right\} \\
        & \dps + \left\{ \int_\om \left[ 2 \frac{\ptl_y V}{V^3}
        \tna V \cdot \na^2_{\tna V} \, I_\dt (\tna V) (\tna V_y)
        (\vy)^2 - 2 \frac{\tna V}{V^2} \na^2_{\tna V} \, I_\dt (\tna V) (\tna V_y) (\vy) \right.
        \right.\\
        & \dps \qquad \qquad + \left. \left. b(\vy)^2 \tna V \cdot \na^2_{\tna V} \, I_\dt (\tna
        V) (\tna V_y) \right] dx\,dy \right\} \\
        & \dps + \left\{ \int_\om
        \tna C \frac{u^\var}{u^\var+\var} V \ptl_{\tna V} \, I_\dt
        (\tna V) (\vy)^2 dx\,dy \right\} \\
        & \dps + \left\{ - \int_\om (\tna b) V_y \, \ptl_{\tna V}
        \, I_\dt (\tna V) (\vy)^2 dx\,dy \right\} \\
        \equiv & \dps \sum^7_{i=1} I_i \, .
        \end{array}
        \end{eqnarray}

        \vskip 5mm

        We now estimate each $I_i$, $(1 \leq i \leq 7)$
        respectively as follows:

        First, due to the facts that $\tilde{\na} {V} |_{y=1} = 0$
        and $\xi$ and $\ptl_\xi I_\dt \geq 0$, one has
        \begin{eqnarray}
        \begin{array}{lcl}
        I_1 & \equiv & \dps \int^L_0 \left[ \ptl_y \left( \frac{\tilde{\na}{V}}{V^2} \right) \ptl_{\tna
        V} \, I_\dt (\tna V) (\vy)^2 - b(\vy)^2 \, \tna V \, \ptl_{\tna V} \, I_\dt (\tna V)  \right. \\
        & & \dps \qquad \qquad + \left. \left. \frac{\tna V}{V^2} \, \ptl_{\tna V} \, I_\dt (\tna V) \, 2
        (\vy) \right] \right|^1_0 dx \\
        & = & \dps - \int^L_0 \left[ \ptl_y \left( \frac{\tna V}{V^2} \right) \ptl_{\tna
        V}\, I_\dt (\tna V) (1+\var)^2 - b(1+\var)^2 \, \tna V \, \ptl_{\tna V} \, I_\dt (\tna
        V) \right. \\
        & & \dps \qquad \qquad + \left. \frac{\tna V}{V^2} \, \ptl_{\tna V} \, I_\dt (\tna
        V)\, 2 (1+\var) \right] (x,0,t) dx\,dt \\
        & \leq & \dps \int^L_0 (1+\var)^2 \left[ \ptl_y \left( \frac{\tna V}{V^2} \right) \ptl_{\tna
        V}\, I_\dt (\tna V) - b \tna V \, \ptl_{\tna V} \, I_\dt (\tna V)
        \right] (x,0,t) dx \\
        & = & \dps + (1+\var)^2 \int^L_0  \left[ \left( \tna V_0 + \tna \left( \frac{\ptl_x P}{U} \right) V \right)
        \cdot \ptl_{\tna V}\, I_\dt (\tna V) \right] \ (x,0,t) dx \\
        & \leq & \tilde{C} \, ,
        \end{array}
        \end{eqnarray}
        here and from now on $\tC = \tC(Q_T)$ is a generic
        positive constant independent $\var$ and $\dt$.  Note that
        in the derivation of (2.34), one has used Lemma 2.2, the
        boundary condition at $y=0$ in (2.9), and the fact that
        $ \dps b(x,0,t)=-\frac{\ptl_x P (x,t)}{U(x,t)}$.  Next, at
        $x=0$, the equation and boundary condition imply that
        $$ \dps ((a+\var)|\ptl_x V|)(x_1 =0,y,t) \leq \frac{|\ptl_t u^\var_1|}{(u^\var_1 + \var)^2} + |\ptl^2_y \, u^\var_1|
        + b(0,y,t) \frac{\ptl_y \, u^\var_1}{(u^\var_1 + \var)^2} + |C| \frac{1}{(u^\var_1 +
        \var)}.$$  Thus, it follows that
        \begin{eqnarray}
        \begin{array}{lcl}
        -I_2 & \equiv & \dps \left.  \int^1_0 (a+\var) I_\dt (\tna V) (\vy)^2 \right|^L_0 dy \\
        & \leq & \dps \int^1_0 \left[ (a+\var) I_\dt (\tna V)
        (\vy)^2 \right] (0,y,t) dy \\
        & \leq & \tC \, ,
        \end{array}
        \end{eqnarray}
        where we have used the assumption (2.3) and the regularity
        assumption on $u^\var_1$.  Next, (2.29) implies that
        \begin{eqnarray}
        \dps I_3 \equiv - \int_\om \frac{\tna V_y}{V^2} \na^2_{\tna V} \,
        I_\dt (\tna V) (\tna V_y)(\vy)^2 dx\,dy \leq 0 \, .
        \end{eqnarray}

        Next, it follows from the structures of $a$, $b$, and $c$
        and Lemma 2.2 that
        \begin{eqnarray}
        \begin{array}{lcl}
        I_4 & \equiv & \dps \int_\om \left[ 2 \left( \frac{\tna V}{V^2} \right) \ptl_{\tna V} \,
        I_\dt (\tna V) - V_x (\tna a) \ptl_{\tna V} \, I_\dt (\tna
        V) (\vy)^2 \right. \\
        & & \dps \qquad + \ptl_x a \, I_\dt (\tna V) (\vy)^2 + \ptl_y (b(\vy)^2) \tna V \, \ptl_{\tna
        V} \, I_\dt (\tna V) \\
        & & \dps \qquad + \left. \{\frac{C u^\var}{u^\var + \var} \tna
        V \cdot \ptl_{\tna V} \, I_\dt (\tna V) - C \var V
        \, \tna V \cdot \ptl_{\tna V} \, I_\dt (\tna V)\} (\vy)^2
        \right] dx\,dy \\
        & \leq & \dps \tC \int_\om (\vy)^2 I_\dt (\tna V) dx\,dy
        \, .
        \end{array}
        \end{eqnarray}

        On the other hand, (2.30) implies that
        \begin{eqnarray}
        \begin{array}{ll}
        I_5 & \equiv \dps \int_\om \tna V \cdot \na^2_{\tna V} \,
        I_\dt (\tna V) (\tna V_y) \left[\{ 2 \frac{\ptl_y V}{V^3}
         - \frac{2}{V^2(\vy)} + b \}(\vy)^2 \right] dx\,dy
        \\
        & \lra 0 \qquad {\rm as} \qquad \dt \lra 0^+ \, .
        \end{array}
        \end{eqnarray}

        Next, it follows from Lemma 2.2 that
        \begin{eqnarray}
        I_6 \equiv \dps \int_\om \tna C \frac{u^\var}{u^\var
        +\var} V \, \ptl_{\tna V} \, I_\dt (\tna V) (\vy)^2 dx\,dy
        \leq \tC \, .
        \end{eqnarray}

        Finally, one has
        \begin{eqnarray}
        \begin{array}{lcl}
        I_7 & \equiv & \dps - \int_\om (\tna b) V_y \, \ptl_{\tna V} \, I_\dt (\tna V) (\vy)^2 dx\,dy
        \\
        & \leq & \dps \tC \int_\om |V_y| (\vy)^2 dx\,dy \\
        & \leq & \dps \tC \int_\om |\ptl_y u^\var (x,y,t)| dx\,dy \, .
        \end{array}
        \end{eqnarray}

        Taking limit as $\dt \ra 0^+$ in (2.33) and using (2.34) -
        (2.40), we get
        \begin{eqnarray*}
        \begin{array}{ll}
        & \dps \frac{d}{dt} \int_\om |\tna V (x,y,t)| (\vy)^2 dx\,dy
        \\
        \leq & \dps \tC + \tC \int_\om |\tna V (x,y,t)| (\vy)^2
        dx\,dy + \tC \int_\om |\ptl_y u^\var (x,y,t)| dx\,dy
        \end{array}
        \end{eqnarray*}
        which yields immediately the desired estimate (2.26) due
        to Lemma 2.4 and Lemma 2.2.  Thus Lemma 2.5
        is proved. \qquad \qquad \qquad \qquad \qquad \qquad \qquad
        \qquad \qquad \qquad \qquad \qquad \qquad \qquad $\square$

    \end{enumerate}
\end{enumerate}

As an immediate consequence of Lemma 2.4 and Lemma 2.5, we have
proved Proposition 2.3.  An immediate corollary of Proposition 2.3
and its proof is that

\vskip 1cm

\begin{corollary}
{\rm For any $\alpha > 0$, there exists a positive constant
$C_5=C_5(\om_T, \alpha)$ independent of $\var$ such that
\begin{eqnarray}
\int\!\!\!\int_{Q_T} (1-y)^\alpha |u^\var_{yy}| dx\,dy\,dt \leq
C_5
\end{eqnarray}
}
\end{corollary}

\vskip 1cm

\begin{enumerate}
    \item[]
    \begin{enumerate}
        \item[{\bf Proof:}] In the proof of Lemma 2.5, if using the
        weight $(\vy)^\alpha$ instead of $(\vy)^2$, one can derive
        similarly that
        \begin{eqnarray}
        \sup_{0 \leq t \leq T} \int_\om (\vy)^\alpha |\tna V
        (x,y,t)| dx\,dy \leq C
        \end{eqnarray}
        which, together with the equation in (2.9), implies
        (2.41).  The proof of Corollary 2.6 is complete.
        \begin{flushright}
        $\square$
        \end{flushright}
    \end{enumerate}
\end{enumerate}

Now we are in the position to prove Theorem 2.1.

\vskip 1cm

\begin{enumerate}
    \item[]
    \begin{enumerate}
        \item[{\bf Proof}] {\bf of Theorem 2.1:} It follows from Lemma
        2.2, Proposition 2.3, Lemma 2.4 and Lemma 2.5 that there
        exist a subsequence of $\{u^\var\}$, which solves (2.9) and
        will be still denoted by $u^\var$, and a function $u
        \in L^\infty (Q_T) \cap BV (Q_T)$ such that as $\var \ra
        0^+$ and $\alpha > 0$,
        \begin{eqnarray}
        u^\var \ra u \quad {\rm a.e.\ in} \quad Q_T, \quad
        (1-y)^{\frac{\alpha -1}{2}} \, u^\var_y \rightharpoonup (1-y)^{\frac{\alpha -1}{2}} \,
        u_y \quad {\rm weakly\ in} \quad Q_T \, ,
        \end{eqnarray}
        and $(1-y)^\alpha \, u^\var_{yy} \ra (1-y)^\alpha \,
        u_{yy}$ is measure on $Q_T$.  Furthermore, the limiting
        function $u(x,y,t)$ satisfies the requirements (2.4),
        (2.5) and (2.6).  Next, we show that $u(x,y,t)$ satisfies
        the equality (2.8).  For any $\vp \in C^1 (\bar{Q}_T)$
        with $\vp|_{t=0} = \vp|_{x=0} = \vp|_{x=L} = 0$ and
        $\alpha > 0$, one has from (2.9) that
        \begin{eqnarray}
        \begin{array}{ll}
        & \dps - \int_\om (u^\var + \var)^{-1} (1-y)^\alpha \, \vp \big|^t
        dx\,dy + \int_{Q_t} (u^\var + \var)^{-1} (1-y)^\alpha \,
        \ptl_t \, \vp \, dx\,dy\,dt \\
        & \dps + \int_{Q_t} (u^\var + \var)^{-1} \, \ptl_x ((1-y)^\alpha (a+\var) \vp)
        dx\,dy\,dt + \int_{Q_t} ((1-y)^\alpha \vp)_y \, u^\var_y \, dx\,dy\,dt \\
        & \dps + \int_{Q_t} ((1-y)^\alpha \, b \vp)_y \, (u^\var + \var)^{-1} dx\,dy\,dt
        + \int^t_0\!\!\!\int^L_0 \left[ \vp \left( u^\var_y + \frac{b}{u^\var + \var} \right) \right]
        (x,0,t) \, dx\,dt \\
        & \dps + \int_{Q_t} (1-y)^\alpha \, \vp c u^\var (u^\var
        + \var)^{-2} dx\,dy\,dt \\
        = & 0 \, .
        \end{array}
        \end{eqnarray}

        Note that $$[u^\var_y + b(u^\var + \var)^{-1}] (x,y=0,t) =
        \ups_0 (x,t)$$ due to the boundary condition at $y=0$ in
        (2.9) and the fact $$b(x,y=0,t) = -(\ptl_x P (x,t))
        U(x,t).$$  Consequently, (2.44) becomes
        \begin{eqnarray*}
        \begin{array}{ll}
        & \dps - \int_\om (u^\var + \var)^{-1} (1-y)^\alpha \, \vp \big|^t
        dx\,dy + \int_{Q_t} (u^\var + \var)^{-1} (1-y)^\alpha \,
        \ptl_t \, \vp \, dx\,dy\,dt \\
        & \dps + \int\!\!\!\int_{Q_t} [(1-y)^\alpha (u^\var + \var)^{-1} \, \ptl_x ((a+\var) \vp)
         + \left((1-y)^\alpha \, \vp_y - \alpha (1-y)^{\alpha -1} \, \vp \right)
        u^\var_y ]\, dx\,dy\,dt \\
        & \dps + \int^t_0\!\!\!\int^L_0 \ups_0 (x,t) \vp (x,0,t) \, dx\,dt
         + \int\!\!\!\int_{Q_t} (1-y)^\alpha \, \vp c u^\var (u^\var + \var)^{-2} dx\,dy\,dt \\
        = & 0 \, .
        \end{array}
        \end{eqnarray*}

        Now, (2.8) follows from passing to the limit $\var \ra
        0^+$ in the above identity and using the convergence
        in (2.43).  Now we verify the boundary conditions.
        First, it follows from $u \in BV(Q_T)$ that $\g u|_{y=1}$
        exists.  Then (2.5) for $u$ implies $\g u|_{y=1} = 0$.
        Next, we consider the boundary condition at $y=0$.
        Observe that $\g u_y |_{y=0}$ exists due to (2.6).

        Define
        $\te \in C^\infty_c [0,1)$ with the properties that $0
        \leq \te \leq 1$, $\te(0)=1$, $\te(y) \equiv 0$ for $y
        \in [\frac{1}{2}, 1)$, and $\int^1_0 \te (y) dy=1$.  Let
        $\psi (x,t)$ be any function in $C^1 ([0,L] \times [0,T])$
        such that $\psi(0,t)=\psi(L,t)=\psi(x,t=0)=0$.  Set
        \begin{eqnarray}
        \vp (x,y,t)= \eta_\dt (y) \, \psi(x,t) \, , \quad {\rm for}
        \quad 0<\dt < 1 \, ,
        \end{eqnarray}
        where
        \begin{eqnarray*}
        \eta_\dt (y) = -\int^1_y \te_\dt (s) ds \, , \quad \te_\dt
        (y) = \dt^{-1} \, \te \left( \frac{y}{\dt} \right) \, .
        \end{eqnarray*}

        Clearly,
        \begin{eqnarray}
        \eta^{'}_\dt (y) = \te_\dt (y) \geq 0 , \quad \eta_\dt (0)=-1 , \quad
        \eta_\dt(y)=0 \quad {\rm for} \quad y \geq
        \frac{1}{2} \dt , \quad {\rm and} \quad |\eta_\dt (y)| \leq 1
        \, .
        \end{eqnarray}

        Using $\vp$ in (2.45) as a test function in (2.8), one
        gets
        \begin{eqnarray*}
        \begin{array}{ll}
        0 = & \dps -\int_\om u^{-1} \vp (1-y)^\alpha dx\,dy \big|^t +
        \int\!\!\!\int_{Q_t} u^{-1} (1-y)^\alpha \, \eta_\dt (y) \, \ptl_t \,
        \psi \, dx\,dy\,dt \\
        & \dps + \int\!\!\!\int_{Q_t} (-\alpha) (1-y)^{\alpha -1} \, \vp \, u_y \, dx\,dy\,dt
        + \int\!\!\!\int_{Q_t} (1-y)^\alpha \, \te_\dt (y) \,
        \psi u_y \, dx\,dy\,dt \\
        & \dps + \int\!\!\!\int_{Q_t} (a \vp)_x \, (1-y)^\alpha \, u^{-1} \,
        dx\,dy\,dt + \int\!\!\!\int_{Q_t} (1-y)^\alpha \, b \, \te_\dt \, \psi \, u^{-1} \,
        dx\,dy\,dt \\
        & \dps + \int\!\!\!\int_{Q_t} ((1-y)^\alpha \, b)_y \, \vp u^{-1} \, dx\,dy\,dt
        + \int^t_0\!\!\int^L_0 \ups_0 (x,t) \, \vp (x,0,t) dx\,dt \\
        & \dps + \int\!\!\!\int_{Q_t} (1-y)^\alpha \, \vp c u^{-1} \,
        dx\,dy\,dt \, .
        \end{array}
        \end{eqnarray*}

        Taking limit as $\dt \ra 0^+$ in the above identity and
        taking into account of (2.46) and Lebesgue dominant
        convergence theorem, one gets that
        $$
        0= \int^t_0\!\!\!\int^L_0 \left\{ \g \, u_y \big|_{y=0}
        (x,t) + \left[ \left( \frac{-\ptl_x P}{U} \right) \left(
        \frac{1}{\g u|_{y=0}} \right) \right] (x,t) - U_0 (x,t)
        \right\} \psi (x,t) dx\,dt
        $$

        It follows that
        $$
        \g \, u_y \big|_{y=0} = \ups_0 + \frac{\ptl_x P}{U}
        \frac{1}{\g u|_{y=0}}
        $$
        which is the desired boundary condition at $y=0$.
        Finally, we verify that $\g u|_{x=0} (y,t) = u_1 (y,t)$,
        a.e..  Indeed, let $\psi (y,t) \in C^1([0,1] \times
        [0,T])$ and $\eta(x) \in C^1 [0,L]$ be arbitrary such
        that $\eta(0)=1$ and $\eta(L)=0$.  Furthermore, set
        $\vp(x,y,t) = \eta(x) \, \psi(y,t)$.  Then one has
        \begin{eqnarray*}
        \begin{array}{ll}
        & \dps \int\!\!\!\int_{Q_T} \vp (x,y,t) \, \ptl_x u (x,y,t)
        dx\,dy\,dt \\
        = & \dps \lim_{\var \ra 0^+} \int\!\!\!\int_{Q_T} \vp
        (x,y,t) \, \ptl_x u^\var (x,y,t) dx\,dy\,dt \\
        = & \dps \lim_{\var \ra 0^+} \left[ - \int^T_0\!\!\!\int^1_0 u^\var_1 (y,t) \, \psi (y,t) dy\,dt
        - \int\!\!\!\int_{Q_T} \ptl_x \, \vp (x,y,t) \, u^\var (x,y,t) dx\,dy\,dt \right] \\
        = & \dps - \int^T_0\!\!\!\int^1_0 u_1 (y,t) \, \psi (y,t) dy\,dt
        - \int\!\!\!\int_{Q_T} \ptl_x \, \vp (x,y,t) \, u (x,y,t) dx\,dy\,dt \\
        = & \dps - \int^T_0\!\!\!\int^1_0 u_1 (y,t) \, \psi (y,t) dy\,dt
        + \int^T_0\!\!\!\int^1_0 \g u \big|_{x=0} (y,t) \, \psi (y,t) dy\,dt\\
        & \dps+ \int\!\!\!\int_{Q_T} \phi (x,y,t) \, \ptl_x \, u (x,y,t) dx\,dy\,dt
        \end{array}
        \end{eqnarray*}

        Consequently,
        $$
        \int^T_0\!\!\!\int^1_0 \g u \big|_{x=0} (y,t) \, \psi (y,t) dy\,dt
        = \int^T_0\!\!\!\int^1_0 u_1 (y,t) \, \psi (x,t) dy\,dt
        $$
        so the desired conclusion follows.  Thus we have shown
        that $u(x,y,t)$ is indeed a weak solution to the
        initial-boundary value problem (2.1).  This completes the
        proof of Theorem 2.1.
        \begin{flushright}
        $\square$
        \end{flushright}
    \end{enumerate}
\end{enumerate}

In the next section, we will show that such a weak solution to
(2.1) is unique.

\vskip 2cm

\section{Uniqueness And Continuous Dependence Of The Weak \\
Solution} \indent

\setcounter{equation}{0}

In this section, we will study the uniqueness and continuous
dependence of weak solutions to the initial-boundary value problem
(1.7) defined in Definition 1.1.  The main strategy is show a
$L^1$-contraction for weak solutions to (1.7) motivated by the
pioneering weak of Krushkov [6].  We start with some basic
structure of a weak solution to (1.7).  Let $u$ be a weak solution
to (1.7).  Then it follows from the definition that $u \in
BV(Q_T) \cap L^\infty (Q_T)$, $\ptl_y u \in L^2(Q_T)$, $\ptl^2_y
\, u \in \mathcal{M}_{loc} (Q_T)$, and hence $\ptl_y u \in
BV^{loc}_y (Q_T)$, here and in what follows, $BV^{loc}_y (Q_T)$
denotes the set of functions $w \in L^1_{loc} (Q_T)$ with $\ptl_y
w$ being a local Random measure on $Q_T$ ([12]).  Let $\G_u$ be the set of approximate jumps of $u$ in $Q_T$
with $\nu=(\nu_x, \nu_y, \nu_t)$ being a normal on $\G_u$.  Denote
by $w^\pm (P)$ be the limiting values of $w$ in the direction $\pm
\nu$ for $P \G_u$.  We also denote by $H$ the 2-dimensional
Hansdorff measure.  Then our first result in this section is on
the jump-conditions for a weak solution $u$.

\setcounter{theorem}{0}

\vskip 1cm

\begin{proposition}
{\rm Let $u$ be a weak solution to the initial-boundary value
problem (1.7).  Then $H$-almost everywhere on $\G_u$, one has that
\begin{eqnarray}
\left( \frac{1}{u^+} - \frac{1}{u^-} \right) (\nu_t + a \nu_x) = 0
, \quad \nu_y = 0 \quad {\rm on} \quad \G_u \, .
\end{eqnarray}
}
\end{proposition}

\vskip 1cm

\begin{enumerate}
    \item[]
    \begin{enumerate}
        \item[{\bf Proof:}] Set $\ups = \frac{1}{u}$.  Then it
        follows from (1.13) in the Definition 1.1 that
        \begin{eqnarray}
        \ptl_t \ups + \ptl^2_y \, u + \ptl_x (a \ups) + \ptl_y (b
        \ups) - (\ptl_x \, a + \ptl_y \, b + c) \ups = 0
        \end{eqnarray}
        in the sense of distribution.  Since $\ups \in BV^{loc}
        (Q_T)$.  Thus (3.2) is in fact a measure equality.  Let
        $S$ be any $H$-measurable subset of $\G_u$.  Then $S$ is
        $\na \ups$-measurable.  Then it follows from (3.2) that
        $S$ is $\ptl^2_y \, u$-measurable.  Furthermore, (3.2)
        implies that
        \begin{eqnarray}
        \int_S (\ups^+ - \ups^-) \nu_t \, dH + \int_S (\ups^+ - \ups^-) a \,
        \nu_x \, dH + \int_S (\ups^+ - \ups^-) b \, \nu_y \, dH = \int_S
        \ptl^2_y \, u \, .
        \end{eqnarray}

        For any given subset $E$ such that $\bar{E} \subset
        \om_T$, we denote by
        $$
        E^{(x,t)} = \{ \, y \eps (0,1) ; \quad {\rm such\ that\ } \,
        (x,y,t) \eps E \, \} \, .
        $$

        Then one can obtain
        \begin{eqnarray*}
        \begin{array}{lcl}
        \dps \int_S u_{yy} & = & \dps \int_{Q_T} \chi_S \, u_{yy} =
        \int^T_0\!\!\!\int^L_0 dx\,dt \int^1_0 \chi_S (x, \cdot, t)
        \, u_{yy} (x, \cdot, t) \\
        & = & \dps \int^T_0\!\!\!\int^L_0 dx\,dt \int_{S^{(x,t)}} u_{yy} (x, \cdot,
        t) = \int^T_0\!\!\!\int^L_0 \sum_{y \eps S^{(x,t)}} \left( u^r_y (x, y,
        t)- u^l_y (x, y, t) \right) dx\,dt
        \end{array}
        \end{eqnarray*}
        where $\chi_S$ is the characteristic function of $S$,
        $u^r_y$ and $u^l_y$ denote the right and left approximate
        limits of $u_y (x,\cdot,t)$ respectively.  Now arguing in a similar way as in the
        proof of Theorem 3.4.1 in [12] (pp. 308 - 309), one obtains
        \begin{eqnarray}
        \int_S u_{yy} = \int_S \left( u^r_y (x,y,t) - u^l_y
        (x,y,t) \right) \, \nu_y \, dH \, .
        \end{eqnarray}

        On the other hand, since $\ptl_y \, u \eps L^2 (Q_T)$, so
        \begin{eqnarray}
        0=\int_S \ptl_y \, u = \int_S (u^+ - u^-) \nu_y \, dH \, .
        \end{eqnarray}

        Since $S$ is an arbitrary subset of $\G_u$, it hence
        follows from (3.5) that $\nu_y = 0$ $H$-almost everywhere
        on $\G_u$.  This together with (3.4) shows that
        $$
        \int_S u_{yy} = 0 \, .
        $$

        Consequently,
        $$
        \int_S (\frac {1}{u^+} - \frac {1}{u^-}) (\nu_t + a \nu_x) \, dH = 0
        $$
        which implies (3.1) since $S$ is an arbitrary subset of
        $\G_u$.  Hence the proof of Proposition 3.1 is completed.
        \begin{flushright}
        $\square$
        \end{flushright}
    \end{enumerate}
\end{enumerate}

We remark here that as an immediate consequence, one can show that
for a weak solution $u$, to (1.7), $\overline{\G}_u$ cannot
contain a 2-dimensional open set.  However, we do not use this
remark later.

Motivated by the arguments of Krushkov [6], we analyze
the ``entropy" $|\ups - k|$ in order to prove the uniqueness and continuous dependence of weak solutions to (1.7).  We define a $C^1$-approximation of
the sign function as for any $\dt > 0$,
\begin{eqnarray}
S_\dt (\te) = \left\{
\begin{array}{cc}
1 , & \te > \dt , \\
\dps \frac{2 \te}{\dt} - \frac{\te^2}{\dt^2} , & 0 \leq \te \leq
\dt , \\
\dps \frac{2 \te}{\dt} + \frac{\te^2}{\dt^2} , & -\dt \leq \te < 0 , \\
-1 , & \te < -\dt ,
\end{array}
\right.
\end{eqnarray}
and $I_\dt (\te) = \int^\te_0 S_\dt (\xi) d \xi$.  Clearly $S_\dt
\, \in C^1 (\mathbb{R}^1)$ and $I_\dt \, \in C^2
(\mathbb{R}^1)$.

Now, for any given functions of $f(s)$ and $w(x,y,t)$, $\hat{f}(\ups)$
will denote the functional superposition defined by Volport in [12], and $\bar{w} (x,y,t)$ will denote the symmetric mean of
$w(x,\cdot ,t)$ as in [12].  For any given constant $k$,
$\hat{I^{'}_\dt} (\ups - k)$ is bounded and measurable with
respect to $\na \ups$, and consequently bounded and measurable
with respect to $u_{yy}$ due to (3.2).  Thus for any smooth test
function $\vp \in C^2_0 (Q_T)$, one has from (3.2) and the
$BV$-calculus in [12] that
\begin{eqnarray}
\begin{array}{ll}
& \dps \int_{Q_T} \left( \phi \{ \ptl_t \, I_\dt (\ups - k) + a \,
\ptl_x \, I_\dt (\ups - k) + b \, \ptl_y \, I_\dt (\ups - k) \} -
c \, \phi \, \widehat{I^{'}_\dt} (\ups - k) \ups \right) \\
= & \dps \int_{Q_T} \left( -\phi \widehat{I^{'}_\dt} (\ups - k) \,
\ptl^2_y \, u \right)
\end{array}
\end{eqnarray}

Note that the singular set of $\ups$ has $H$-measure zero, and on
the set of approximate continuous points of $\ups$, $\widehat{I^{'}_\dt} (\ups - k)
= \overline{\widehat{I^{'}_\dt}} (\ups - k) = I^{'}_\dt
(\overline{\ups} - k)$. It follows from Proposition 3.1, (3.6),
and the $BV$-calculus that
\begin{eqnarray}
\begin{array}{lcl}
\dps \int_{Q_T} \left( -\phi \, \widehat{I^{'}_\dt} (\ups - k) \,
\ptl^2_y \, u \right) & = & \dps \int_{Q_T} \left( -\phi \,
\overline{\widehat{I^{'}_\dt}} (\ups - k) \, \ptl^2_y \, u \right)
+ \int_{Q_T} -\phi \left( \widehat{I^{'}_\dt} (\ups - k) -
\overline{\widehat{I^{'}_\dt}} (\ups - k)
\right) \, \ptl^2_y \, u \\
& = & \dps \int_{Q_T} (-\phi) \, \ptl_y \left(\widehat{I^{'}_\dt}
(\ups - k) \, \ptl_y \, u \right) + \int_{Q_T} \phi \,
\widehat{\widehat{I^{''}_\dt}}
(\ups - k) \, \ptl_y \, \ups \, \overline{\ptl_y \, u} \, dx\,dy \\
& = & \dps \int_{Q_T} \ptl_y \, \phi \, \widehat{I^{'}_\dt} (\ups
- k) \, \ptl_y \, u \, dx\,dy\,dt + \int_{Q_T} \phi \, I^{''}_\dt
(\ups - k) \, \ptl_y \, \ups \, \ptl_y \, u \, dx\,dy\,dt
\end{array}
\end{eqnarray}
where we have used the fact that $\ptl_y \, u \in L^2(Q_T)$ and
$\ptl_y \, \ups \in L^2_{loc} (Q_T)$, which follows from the
definition.  Thus, it follows from (3.7), (3.8) and Proposition
3.1 again that
\begin{eqnarray}
\begin{array}{ll}
& \dps \int_{Q_T} \left\{ I_\dt (\ups - k) \left( \ptl_t \, \phi +
\ptl_x (a \phi) + \ptl_y (b \phi) \right) + I^{'}_\dt (\ups - k)
\, \ups c \phi \right\} dx\,dy\,dt \\
= & \dps - \int_{Q_T} \ptl_y \, \phi \, \widehat{I^{'}_\dt} (\ups
- k) \, \ptl_y \, u \, dx\,dy\,dt - \int_{Q_T} \phi \, I^{''}_\dt
(\ups - k) \, \ptl_y \, \ups \, \ptl_y \, u \, dx\,dy\,dt
\end{array}
\end{eqnarray}

\vskip 5mm

Consequently, we have shown the following lemma.

\vskip 1cm

\begin{lemma}
{\rm Let $u$ be a weak solution to (2.1) and $\ups = u^{-1}$. Then
for any constant $k$ and $\phi \in C^\infty_0 (Q_T)$, one has the
identity
\begin{eqnarray}
\begin{array}{ll}
& \dps \int_{Q_T} \left\{ |\ups - k| \left( \ptl_t \, \phi +
\ptl_x (a \phi) + \ptl_y (b \phi) \right) + {\rm sgn\ } (\ups - k)
\, \ups c \phi \right\} dx\,dy\,dt \\
= & \dps - \int_{Q_T} (\ptl_y \, \phi) \, {\rm sgn\ } (\ups - k)
\, \ptl_y \, u \, dx\,dy\,dt - \lim_{\dt \ra 0^+} \int_{Q_T} \phi
\, I^{''}_\dt (\ups - k) \ptl_y \, \ups \, \ptl_y \, u \,
dx\,dy\,dt
\end{array}
\end{eqnarray}
}
\end{lemma}

\vskip 5mm

As an immediate consequence of this lemma, one can have

\vskip 1cm

\begin{corollary}
{\rm Let $u$ be a weak solution to (2.1) and $\ups=u^{-1}$.  Then
for any constant $k$, and $\psi \in C^\infty_0 (Q_T)$, $\psi \geq
0$, it holds that
\begin{eqnarray}
\int_{Q_T} \{ \, |\ups - k| (\ptl_t \, \psi + \ptl_x \, (a \psi) +
\ptl_y (b \psi)) + \ptl_y \, \psi \, {\rm sgn} (\ups - k) \ptl_y
\, u - \, {\rm sgn} (\ups - k) \ups c \psi \, \} dx\,dy\,dt \geq 0
\, .
\end{eqnarray}

Now let $u_1$ and $u_2$ be two weak solutions to (2.1) with their
corresponding initial and boundary data respectively.  Set $\ups_i
= (u_i)^{-1}$.  Then our key estimate in the argument for
continuous dependence is the following proposition: }
\end{corollary}

\vskip 1cm

\begin{proposition}
{\rm Let $u_1$ and $u_2$ be two weak solutions to (2.1). For any $\phi \in C^\infty_0 (Q_T)$, it holds that
\begin{eqnarray}
\int_{Q_T} \{ \, |\ups_1 - \ups_2| (\ptl_t \, \phi + \ptl_x \, (a
\phi) + \ptl_y (b \phi) + c \phi) - \ptl_y \, \phi \, \ptl_y |u_1
- u_2| \, \} dx\,dy\,dt \geq 0 \, .
\end{eqnarray}
}
\end{proposition}

\vskip 1cm

\begin{enumerate}
    \item[]
    \begin{enumerate}
        \item[{\bf Proof:}] This proposition will be proved by
        modifying the doubling-variable argument of Krushkov [6]
        based on Lemma 3.2.  Denote by $\rx = (x,y,t)$ and
        $\brx = (\bar{x}, \bar{y},
        \bar{t})$. Now for any nonnegative function $\psi =
        \psi (\rx, \brx) \in C^\infty$
        such that
        \begin{eqnarray}
        \mathrm{supp} \Psi (\cdot, \brx) \subset\subset Q_T, \quad
        \forall \, \brx \, \in \, Q_T, \quad \mathrm{supp} \Psi (\rx, \cdot)
        \subset\subset Q_T, \quad \forall \, \rx \, \in \, Q_T \, ,
        \end{eqnarray}
        it follows from (3.10) in Lemma 3.2 that
        \begin{eqnarray}
        \begin{array}{ll}
        & \dps \int_{Q_T}\!\!\!\int_{Q_T} \{ \, |\ups_1(X) - \ups_2(\bar{X})| \left(\ptl_t \, \Psi + \ptl_x \,
        (a(X) \Psi) + \ptl_y (b(X) \Psi) + c(X) \, \ups_1(X) \right) \, \mathrm{sgn} (\ups_1(X) - \ups_2(\bar{X}))
        \Psi \\
        & \qquad \qquad \dps + \ptl_y \, \Psi \, \mathrm{sgn} (\ups_1(X) -
        \ups_2(\bar{X})) \, \ptl_y \, u_1(X) \, \} dX\,d\bar{X} \\
        = & \dps - \lim_{\dt \ra 0^+} \int_{Q_T}\!\!\!\int_{Q_T} \Psi I^{''}_\dt (\ups_1(X) -
        \ups_2(\bar{X})) \, \ptl_y \, \ups_1 (X) \, \ptl_y \, u_1(X) \, dX\,d\bar{X} \, ,
        \end{array}
        \end{eqnarray}
        and
        \begin{eqnarray}
        \begin{array}{ll}
        & \dps \int_{Q_T}\!\!\!\int_{Q_T} \{ \, |\ups_2(\bar{X}) - \ups_1(X)| \big(\ptl_{\bar{t}} \, \Psi
        + \ptl_{\bar{x}} \, (a(\bar{X}) \Psi) + \ptl_{\bar{y}} (b(\bar{X}) \Psi) \big)
        + c(\bar{X}) \, \ups_2(\bar{X}) \, \mathrm{sgn} (\ups_2(\bar{X}) -
        \ups_1(\bar{X})) \Psi \\
        & \qquad \qquad \dps + \ptl_{\bar{y}} \, \Psi \, \mathrm{sgn} (\ups_2(\bar{X}) - \ups_1(X))
        \ptl_{\bar{y}} \, u_2(\bar{X}) \, \} d\bar{X}\,dX \\
        = & \dps - \lim_{\dt \ra 0^+} \int_{Q_T}\!\!\!\int_{Q_T} \Psi I^{''}_\dt (\ups_2(\bar{X}) -
        \ups_1(X)) \ptl_{\bar{y}} \, \ups_2 (\bar{X}) \,
        \ptl_{\bar{y}} \, u_2(\bar{X}) \, d\bar{X}\,dX \, ,
        \end{array}
        \end{eqnarray}

        Adding (3.14) to (3.15) shows that
        \begin{eqnarray}
        \begin{array}{ll}
        & \dps \int_{Q_T}\!\!\!\int_{Q_T} \{|\ups_1(X) - \ups_2(\bar{X})|
        \big( (\ptl_t \Psi + \ptl_{\bar{t}} \Psi) + \ptl_x (a(X) \Psi)
        + \ptl_{\bar{x}} (a(\bar{X}) \Psi) + \ptl_y (b(X) \Psi) +
        \ptl_{\bar{y}} (b(\bar{X}) \Psi) \\
        & \qquad \qquad \qquad \qquad \qquad \qquad \dps + c(X) \Psi \big) -
        ( c(X) - c(\bar{X})) \, \Psi \ups_2(\bar{X}) \, \mathrm{sgn} (\ups_2(\bar{X}) -
        \ups_1(X)) \} dX\,d\bar{X} \\
        & \dps + \int_{Q_T}\!\!\!\int_{Q_T} \{ \, \ptl_y \, \Psi \, \mathrm{sgn} (\ups_1(X) -
        \ups_2(\bar{X})) \ptl_y \, u_1(X) + \ptl_{\bar{y}} \, \Psi \, \mathrm{sgn} (\ups_2(\bar{X}) -
        \ups_1(X)) \ptl_{\bar{y}} \, u_2(\bar{X}) \} dX\,d\bar{X} \\
        = & \dps - \lim_{\dt \ra 0^+} \int_{Q_T}\!\!\!\int_{Q_T} \Psi I^{''}_\dt (\ups_2(\bar{X}) - \ups_1(X))
        (\ptl_y \, \ups_1(X) \, \ptl_y \, u_1(X) + \ptl_{\bar{y}} \, \ups_2(\bar{X})
        \, \ptl_{\bar{y}} \, u_2(\bar{X}) \, dX\,d\bar{X} \, ,
        \end{array}
        \end{eqnarray}
        where we have used the property that $I^{''}_\dt (\xi)$ is
        an even function.  We will estimate each term in (3.16)
        respectively for a special choice of the test function
        $\Psi$.  Let $\rho \in C^\infty_0 (\mathbb{R}^1)$, $\rho
        (\xi) \geq 0$, $\mathrm{supp} \, \rho (\cdot) \subset
        (-1,1)$, and $\int_{\mathbb{R}^1} \rho (\xi) d\xi = 1$.
        For any $h>0$ sufficiently small, set $\rho_h (\xi) =
        \frac{1}{h} \rho (\frac{\xi}{h})$.

        Set
        \begin{eqnarray}
        \Psi(X,\bar{X}) = \frac{1}{8} \, \phi \left(\frac{x+\bar{x}}{2}, \frac{y+\bar{y}}{2}, \frac{t+\bar{t}}{2}
        \right) \rho_h \left(\frac{x-\bar{x}}{2}\right) \rho_h \left(\frac{y-\bar{y}}{2}
        \right) \rho_h \left(\frac{t-\bar{t}}{2} \right), \,
        \forall (X,\bar{X}) \eps Q_T \times Q_T \quad
        \end{eqnarray}
        which clearly satisfies the requirements in (3.13).  In the
        following, for the convenience of notations, we will use
        $\ptl_1$, $\ptl_2$, $\ptl_3$ to denote $\ptl_x$, $\ptl_y$,
        $\ptl_t$ respectively, and $$\tilde{\rho}_h (X,\stackrel{\rightharpoonup}{X}) =
        \frac{1}{8} \rho_h \left( \frac{x-\bar{x}}{2} \right) \, \rho_h \left( \frac{y-\bar{y}}{2}
        \right) \, \rho_h \left( \frac{t-\bar{t}}{2} \right).$$
        Since $$\ptl_t \Psi + \ptl_{\bar{t}} \Psi = \ptl_3 \phi \left(\frac{X+\bar{X}}{2}\right)
        \tilde{\rho}_h (X,\bar{X}),$$ so it follows from the
        Lebesgue dominant convergence theorem as in [6] that
        \begin{eqnarray}
        \int_{Q_T}\!\!\!\int_{Q_T} |\ups_1(X) - \ups_2(\bar{X})|
        \ptl_t(\Psi) + \ptl_{\bar{t}}\Psi) dX\,d\bar{X}
        \lra \int_{Q_T} |\ups_1(X) - \ups_2(X)|
        (\ptl_t \, \phi (X) \, dX , \, h \ra 0^+ \, . \, \quad
        \end{eqnarray}

        Next, one has
        \begin{eqnarray}
        \begin{array}{ll}
        & \dps \int_{Q_T}\!\!\!\int_{Q_T} |\ups_1(X) - \ups_2(\bar{X})|
        (\ptl_x(a(X)\Psi) + \ptl_{\bar{x}}(a(\bar{X})\Psi)
        dX\,d\bar{X} \\
        = & \dps \int_{Q_T}\!\!\!\int_{Q_T} (\ptl_x \, a(X) + \ptl_{\bar{x}} a(\bar{X})
        |\ups_1(X) - \ups_2(\bar{X})|\, \Psi \, dX\,d\bar{X} \\
        & \dps + \int_{Q_T}\!\!\!\int_{Q_T} |\ups_1(X) -
        \ups_2(\bar{X})| \, \frac{1}{2} (a(X)+a(\bar{X})) \,
        \ptl_1 \, \phi \left( \frac{X+\bar{X}}{2} \right) \,
        \stackrel{\rightharpoonup}{\rho_h}\left( \frac{X-\bar{X}}{2} \right)
        \, dX\,d\bar{X} \\
        & \dps + \int_{Q_T}\!\!\!\int_{Q_T} |\ups_1(X) -
        \ups_2(\bar{X})| \, \frac{1}{2} (a(X)-a(\bar{X})) \,
        \left( \frac{\rho^{'}_h}{\rho_h} \right) \,
        \left( \frac{x-\bar{x}}{2} \right) \, \Psi
        \, dX\,d\bar{X} \\
        \equiv & \dps \sum^3_{i=1} \, I_i \, .
        \end{array}
        \end{eqnarray}

        Clearly,
        \begin{eqnarray}
        \begin{array}{lcl}
        I_1 & \equiv & \dps \int_{Q_T}\!\!\!\int_{Q_T} |\ups_1(X) -
        \ups_2(\bar{X})| (\ptl_x \, a(X) + \ptl_{\bar{x}} \, a(\bar{X})) \Psi \, dX\,d\bar{X} \\
        & & \dps \lra \, 2 \int_{Q_T} |\ups_1(X) -
        \ups_2(X)| (\ptl_x \, a(X)) \phi(X) \, dX , \quad
        {\rm as} \quad h \ra 0^+ \, ,
        \end{array}
        \end{eqnarray}
        \begin{eqnarray}
        \begin{array}{lcl}
        I_2 & \equiv & \dps \int_{Q_T}\!\!\!\int_{Q_T} |\ups_1(X) -
        \ups_2(\bar{X})| \, \frac{1}{2} (a(X) + a(\bar{X})) \, \ptl_1 \, \phi
        \left( \frac{X+\bar{X}}{2} \right) \, \tilde{\rho}_h
        \left( \frac{X-\bar{X}}{2} \right) dX\,d\bar{X} \\
        & & \dps \lra \int_{Q_T} |\ups_1(X) -
        \ups_2(X)| a(X) \, \ptl_x \, \phi(X) \, dX , \quad
        {\rm as} \quad h \ra 0^+ \, .
        \end{array}
        \end{eqnarray}

        We rewrite $I_3$ as
        \begin{eqnarray}
        \begin{array}{lcl}
        I_3 & \equiv & \dps \int_{Q_T}\!\!\!\int_{Q_T} |\ups_1(X) -
        \ups_2(\bar{X})| \, \frac{1}{2} (a(X) - a(\bar{X})) \left( \frac{\rho^{'}_h}{\rho_h} \right) \,
        \left( \frac{x-\bar{x}}{2} \right) \, \Psi \,
        dX\,d\bar{X}\\
        & = & \dps \int_{Q_T}\!\!\!\int_{Q_T} \frac{1}{2} \left\{ \, |\ups_1(X) -
        \ups_2(X)| + |\ups_1(\bar{X}) -
        \ups_2(\bar{X})| \right. \\
        & & \dps \quad \quad + \left. \left[ \, (|\ups_1(X) -
        \ups_2(\bar{X})| - |\ups_1(X) - \ups_2(X)| ) + (|\ups_1(X) -
        \ups_2(\bar{X})| - |\ups_1(\bar{X}) - \ups_2(\bar{X})|)
        \right] \right\} \cdot \\
        & & \dps \quad \quad \left[ \frac{1}{2} (a(X) - a(\bar{X})) \right]
        \left( \frac{\rho^{'}_h}{\rho_h} \right) \left( \frac{x-\bar{x}}{2} \right) \, \Psi \,
        dX\,d\bar{X}\\
        & \equiv & \dps \sum^3_{i=1} I_{3i} \, .
        \end{array}
        \end{eqnarray}

        Since,
        \begin{eqnarray*}
        \begin{array}{lcl}
        I_{31} & = & \dps \int_{Q_T}\!\!\!\int_{Q_T} \frac{1}{2} |\ups_1(X) -
        \ups_2(X)| \, \frac{1}{2} (a(X) - a(\bar{X})) \left( \frac{\rho^{'}_h}{\rho_h} \right) \,
        \left( \frac{x-\bar{x}}{2} \right) \, \Psi \, dX\,d\bar{X}\\
        & = & \dps \int_{Q_T}\!\!\!\int_{Q_T} \frac{1}{2} |\ups_1(X) -
        \ups_2(X)| \, \frac{1}{2} (a(X) - a(\bar{X})) \, \ptl_1 \, \phi
        \left(\frac{X+\bar{X}}{2}\right)
        \tilde{\rho}_h \left(\frac{X-\bar{X}}{2} \right) \, dX\,d\bar{X} \\
        & & \dps - \frac{1}{2} \int_{Q_T}\!\!\!\int_{Q_T} |\ups_1(X) -
        \ups_2(X)| \, \ptl_{\bar{x}} \, a(\bar{X}) \Psi \,
        d\bar{X}\,dX \, ,
        \end{array}
        \end{eqnarray*}
        so,
        \begin{eqnarray}
        I_{31} \lra - \frac{1}{2} \int_{Q_T} |\ups_1(X) -
        \ups_2(X)| \, \ptl_{x} \, a(X) \phi(X) \, dX \quad {\rm
        as} \quad h \ra 0^+ \, .
        \end{eqnarray}

        Similarly,
        \begin{eqnarray}
        \begin{array}{lcl}
        I_{32} & \equiv & \dps \int_{Q_T}\!\!\!\int_{Q_T} \frac{1}{2} |\ups_1(\bar{X}) -
        \ups_2(\bar{X})| \, \frac{1}{2} (a(X) - a(\bar{X})) \left( \frac{\rho^{'}_h}{\rho_h} \right) \,
        \left( \frac{x-\bar{x}}{2} \right) \, \Psi \, dX\,d\bar{X}\\
        & & \dps \lra - \frac{1}{2} \int_{Q_T} |\ups_1(X) -
        \ups_2(X)| \, \ptl_{x} \, a(X) \phi(X) \, dX \quad {\rm
        as} \quad h \ra 0^+ \, .
        \end{array}
        \end{eqnarray}

        Note that
        $$\begin{array}{ll}
        |(|\ups_1(X) - \ups_2(\bar{X})| - |\ups_1(X) - \ups_2(X)|) + (|\ups_1(X) -
        \ups_2(\bar{X})| - |\ups_1(\bar{X}) - \ups_2(\bar{X})|)| \\  \dps   \leq |\ups_1(X) -
        \ups_1(\bar{X})| + |\ups_2(X) - \ups_2(\bar{X})|
        \end{array}
        $$ and $a(\cdot)$ is
        smooth.  It thus follows from the Lebesgue dominant
        convergence theorem that as $h \ra 0^+$,
        \begin{eqnarray}
        \begin{array}{lcl}
        I_{33} & \equiv & \dps \int_{Q_T}\!\!\!\int_{Q_T} [ (|\ups_1(X) -
        \ups_2(\bar{X})| - |\ups_1(X) - \ups_2(X)| ) + (|\ups_1(X) -
        \ups_2(\bar{X})| - |\ups_1(\bar{X}) - \ups_2(\bar{X})|) ]
        \\
        & & \dps \quad \quad \frac{1}{2} (a(X) - a(\bar{X})) \left( \frac{\rho^{'}_h}{\rho_h} \right) \,
        \left( \frac{x-\bar{x}}{2} \right) \, \Psi \, dX\,d\bar{X}
        \lra 0 \, .
        \end{array}
        \end{eqnarray}

        Consequently, we have shown from (3.22) - (3.25) that
        $$
        I_3 \lra - \int_{Q_T} |\ups_1(X) -
        \ups_2(X)| \, \ptl_{x} \, a(X) \phi(X) \, dX \quad {\rm
        as} \quad h \ra 0^+ \, .
        $$

        This, together with (3.20) and (3.21), shows that
        \begin{eqnarray}
        \begin{array}{cl}
        & \dps \int_{Q_T}\!\!\!\int_{Q_T} |\ups_1(X) -
        \ups_2(\bar{X})| \, (\ptl_{x} \, (a(X) \Psi) + \ptl_{\bar{x}} \, (a(\bar{X})
        \Psi)) dX\,dX \\
        \lra & \dps \int_{Q_T} |\ups_1(X) -
        \ups_2(X)| \, \ptl_{x} \, (a(X) \phi(X)) \, dX \quad {\rm
        as} \quad h \ra 0^+ \, .
        \end{array}
        \end{eqnarray}

        Similarly, one can show that
        \begin{eqnarray}
        \begin{array}{cl}
        & \dps \int_{Q_T}\!\!\!\int_{Q_T} |\ups_1(X) -
        \ups_2(\bar{X})| \, (\ptl_{y} \, (b(X) \Psi) + \ptl_{\bar{y}} \, (b(\bar{X})
        \Psi)) dX\,d\bar{X} \\
        \lra & \dps \int_{Q_T} |\ups_1(X) -
        \ups_2(X)| \, \ptl_{y} \, (b(X) \phi(X)) \, dX \quad {\rm
        as} \quad h \ra 0^+ \, .
        \end{array}
        \end{eqnarray}

        Next, one has
        \begin{eqnarray}
        \begin{array}{cl}
        & \dps \int_{Q_T}\!\!\!\int_{Q_T} |\ups_1(X) -
        \ups_2(\bar{X})| \, (-c(X) \Psi) \, dX\,dX \\
        \lra & \dps - \int_{Q_T} |\ups_1(X) -
        \ups_2(X)| \, c(X) \phi(X) \, dX \quad {\rm
        as} \quad h \ra 0^+ \, ,
        \end{array}
        \end{eqnarray}
        and
        \begin{eqnarray}
        \begin{array}{cl}
        & \dps \int_{Q_T}\!\!\!\int_{Q_T} |\ups_1(X) -
        \ups_2(\bar{X})| \, (c(\bar{X}) - c(X)) \Psi \, \ups_2(\bar{X}) \, {\rm sgn\ }
         (\ups_2(\bar{X}) - \ups_1(X)) dX\,d\bar{X} \\
        \lra & \dps 0 \qquad \qquad {\rm as} \quad h \ra 0^+ \,
        .
        \end{array}
        \end{eqnarray}

        Finally, we treat the last term on the left hand side of
        (3.16).  Note that
        \begin{eqnarray}
        \begin{array}{cl}
        & \dps \int_{Q_T}\!\!\!\int_{Q_T} \{ \, \ptl_y \, \Psi \, {\rm sgn\ } (\ups_1(X) -
        \ups_2(\bar{X})) \, \ptl_y \, u_1(X) + \ptl_{\bar{y}} \, \Psi \, {\rm sgn\ }
        (\ups_2(\bar{X}) - \ups_1(X)) \, \ptl_{\bar{y}} \, u_2(\bar{X}) \} \,
        dX\,d\bar{X}\\
        = & \dps \lim_{\dt \ra 0^+} \int_{Q_T}\!\!\!\int_{Q_T} \{ \, \ptl_y \, \Psi \, I^{'}_\dt (\ups_1(X) -
        \ups_2(\bar{X})) \, \ptl_y \, u_1(X) + \ptl_{\bar{y}} \, \Psi \,
        I^{'}_\dt (\ups_2(\bar{X}) - \ups_1(X)) \, \ptl_{\bar{y}} \, u_2(\bar{X}) \} \,
        dX\,d\bar{X} \\
        = & \dps \lim_{\dt \ra 0^+} \int_{Q_T}\!\!\!\int_{Q_T} (I^{'}_\dt \, (\ups_1(X) -
        \ups_2(\bar{X})) \, \ptl_y \, u_1(X) + I^{'}_\dt (\ups_2(\bar{X}) - \ups_1(X)) \,
        \ptl_{\bar{y}} \, u_2(\bar{X}) ) \\
        & \dps \qquad \qquad \qquad \qquad \frac{1}{2} \ptl_2 \,
        \phi \left( \frac{X+\bar{X}}{2} \right) \tilde{\rho}_h
        \left( \frac{X-\bar{X}}{2} \right) dX\,d\bar{X} \\
        & \dps + \lim_{\dt \ra 0^+} \int_{Q_T}\!\!\!\int_{Q_T} (I^{'}_\dt \, (\ups_1(X) -
        \ups_2(\bar{X})) \, \ptl_y \, u_1(X) - I^{'}_\dt (\ups_2(\bar{X}) - \ups_1(X)) \,
        \ptl_{\bar{y}} \, u_2(\bar{X}) ) \\
        & \dps \qquad \qquad \qquad \qquad \frac{1}{2} \left( \frac{\rho^{'}_h}{\rho_h} \right)
        \left( \frac{y-\bar{y}}{2} \right) \Psi \, dX\,d\bar{X} \\
        \equiv & \dps J_1 + J_2
        \end{array}
        \end{eqnarray}

        Using the fact that $\ptl_y \, u_1 (\cdot) \in L^2 (Q_T)$
        and $\ptl_{\bar{y}} \, u(\cdot) \in L^2 (Q_T)$, one can
        show easily that
        \begin{eqnarray*}
        \begin{array}{lcl}
        J_1 & = & \dps \lim_{\dt \ra 0^+} \int_{Q_T}\!\!\!\int_{Q_T} (I^{'}_\dt \, (\ups_1(X) -
        \ups_2(\bar{X})) \, \ptl_y \, u_1(X) + I^{'}_\dt (\ups_2(\bar{X}) - \ups_1(X)) \,
        \ptl_{\bar{y}} \, u_2(\bar{X}) ) \\
        & & \dps \qquad \qquad \qquad \qquad \frac{1}{2} \ptl_2 \,
        \phi \left( \frac{X+\bar{X}}{2} \right) \tilde{\rho}_h
        \left( \frac{X-\bar{X}}{2} \right) dX\,d\bar{X} \\
        & = & \dps \int_{Q_T}\!\!\!\int_{Q_T} ( {\rm{sgn\ }} (\ups_1(X) -
        \ups_2(\bar{X})) \, \ptl_y \, u_1(X) + {\rm{sgn\ }} (\ups_2(\bar{X}) - \ups_1(X)) \,
        \ptl_{\bar{y}} \, u_2(\bar{X}) ) \\
        & & \dps \qquad \qquad \qquad \qquad \frac{1}{2} \ptl_2 \,
        \phi \left( \frac{X+\bar{X}}{2} \right) \tilde{\rho}_h
        \left( \frac{X-\bar{X}}{2} \right) dX\,d\bar{X} \\
        & = & \dps \int_{Q_T}\!\!\!\int_{Q_T} (-\ptl_y |u_1(X) -
        u_2(\bar{X})| - \ptl_{\bar{y}} |u_2(\bar{X}) - u_1(X)|) \\
        & & \dps \qquad \qquad \qquad \qquad \frac{1}{2} \ptl_2 \,
        \phi \left( \frac{X+\bar{X}}{2} \right) \tilde{\rho}_h
        \left( \frac{X-\bar{X}}{2} \right) dX\,d\bar{X} \, ,
        \end{array}
        \end{eqnarray*}
        and so,
        \begin{eqnarray}
        J_1 \lra - \frac{1}{2} \int_{Q_T} \ptl_y |
        u_1(X) - u_2(X)| \ptl_y \, \phi(X) \, dX \qquad {\rm as} \quad h
        \ra 0^+ \, .
        \end{eqnarray}

        It remains to estimate $J_2$.  It follows from the fact
        that $\ptl_y \, u_2(\cdot) \in L^2 (Q_T)$, $\ptl_y \, \ups_i(\cdot) \in L^2_{loc}
        (Q_T)$, $\ptl_{\bar{y}} \, u_i(\cdot) \in L^2 (Q_T)$,
        $\ptl_{\bar{y}} \, \ups_i(\cdot) \in L^2 (Q_T)$, and
        Proposition 3.1 that
        \begin{eqnarray}
        \begin{array}{lcl}
        J_2 & = & \dps \lim_{\dt \ra 0^+} \int_{Q_T}\!\!\!\int_{Q_T} (I^{'}_\dt \, (\ups_1(X) -
        \ups_2(\bar{X})) \, \ptl_y \, u_1(X) - I^{'}_\dt (\ups_2(\bar{X}) - \ups_1(X)) \,
        \ptl_{\bar{y}} \, u_2(\bar{X}) ) \\
        & & \dps \qquad \qquad \qquad \qquad \frac{1}{2} \left( \frac{\rho^{'}_h}{\rho_h} \right)
        \left( \frac{y-\bar{y}}{2} \right) \Psi \, dX\,d\bar{X} \\
        & = & \dps \lim_{\dt \ra 0^+} \int_{Q_T}\!\!\!\int_{Q_T} (I^{'}_\dt \, (\ups_1(X) -
        \ups_2(\bar{X})) \, \ptl_y \, u_1(X) + I^{'}_\dt (\ups_2(\bar{X}) - \ups_1(X)) \,
        \ptl_{\bar{y}} \, u_2(\bar{X}) ) \\
        & & \dps \qquad \qquad \qquad \qquad \frac{1}{2} \, \ptl_2 \, \phi \left( \frac{X+\bar{X}}{2} \right)
        \tilde{\rho}_h \left( \frac{X-\bar{X}}{2} \right) dX\,d\bar{X} \\
        & & + \dps \lim_{\dt \ra 0^+} \int_{Q_T}\!\!\!\int_{Q_T} (I^{''}_\dt \, (\ups_1(X) -
        \ups_2(\bar{X})) \, (- \ptl_{\bar{y}} \, u_2(\bar{X})) \ptl_y \, u_1(X) \\
        & & \qquad \qquad \qquad \qquad + I^{''}_\dt
        (\ups_2(\bar{X}) - \ups_1(X)) \,
        (- \ptl_y \, \ups_1(X) ) \ptl_{\bar{y}} \, u_2(\bar{X})) \Psi \, dX\,d\bar{X} \\
        & \equiv & J_1 + J_{22} \, .
        \end{array}
        \end{eqnarray}

        To treat $J_{22}$, we can estimate it by using the facts
        that $\ptl_y \, u_i \, \in L^2 (Q_T)$, $\ptl_{\bar{y}} \, \ups_i \, \in L^2 (Q_T)$,
        $\ptl_y \, \ups_i \, \in L^2_{loc} (Q_T)$, and
        $\ptl_{\bar{y}} \, \ups_i \, \in L^2_{loc} (Q_T)$ for
        $i=1,2$, as
        \begin{eqnarray}
        \begin{array}{lcl}
        J_{22} & = & \dps \lim_{\dt \ra 0^+} \int_{Q_T}\!\!\!\int_{Q_T} (I^{''}_\dt \, (\ups_1(X) -
        \ups_2(\bar{X})) \, (- \ptl_{\bar{y}} \, \ups_2(\bar{X}) \ptl_y \, u_1(X) -
        \ptl_y \, \ups_1(X) \ptl_{\bar{y}} \, u_2(\bar{X})) \Psi
        \, dX\,d\bar{X} \\
        & = & \dps \lim_{\dt \ra 0^+} \int_{Q_T}\!\!\!\int_{Q_T} I^{''}_\dt \, (\ups_1(X) -
        \ups_2(\bar{X})) \, \left( \frac{\ups_2(\bar{X})}{\ups_1(X)} +
        \frac{\ups_1(X)}{\ups_2(\bar{X})} \right) (\ptl_y \, l_n \, u_1(X)) (\ptl_{\bar{y}} \, l_n \,
        u_2(\bar{X}))
        \Psi \, dX\,d\bar{X} \\
        & \leq & \dps \lim_{\dt \ra 0^+} \int_{Q_T}\!\!\!\int_{Q_T} I^{''}_\dt \, (\ups_1(X) -
        \ups_2(\bar{X})) \left( \frac{\ups_2(\bar{X})}{\ups_1(X)} +
        \frac{\ups_1(X)}{\ups_2(\bar{X})} \right) \\
        & & \dps \qquad \qquad \qquad \qquad \frac{1}{2} \left[ (\ptl_y \, l_n \, u_1(X))^2 +
        (\ptl_{\bar{y}} \, l_n \, u_2(\bar{X}))^2 \right]
        \Psi \, dX\,d\bar{X} \, ,
        \end{array}
        \end{eqnarray}

        Finally, we rewrite the integral on the right hand side of
        (3.16) as
        \begin{eqnarray*}
        \begin{array}{lcl}
        & & \dps - \lim_{\dt \ra 0^+} \int_{Q_T}\!\!\!\int_{Q_T} \Psi I^{''}_\dt \, (\ups_2(\bar{X}) -
        \ups_1(X)) \, (\ptl_y \, \ups_1(X) \, \ptl_y \, u_1(X) +
        \ptl_{\bar{y}} \, \ups_2(\bar{X}) \, \ptl_{\bar{y}} \, u_2(\bar{X}))
        \, dX\,d\bar{X} \\
        & = & \dps \int_{Q_T}\!\!\!\int_{Q_T} I^{''}_\dt \, (\ups_2(\bar{X}) -
        \ups_1(X)) \, \left[ (\ptl_y \, l_n \, u_1(X))^2 + (\ptl_{\bar{y}} \, l_n \,
        u_2(\bar{X}))^2 \right] \Psi \, dX\,d\bar{X}
        \Psi \, dX\,d\bar{X} \, ,
        \end{array}
        \end{eqnarray*}

        Thus,
        \begin{eqnarray}
        \begin{array}{lcl}
        & & \dps - \lim_{\dt \ra 0^+} \int_{Q_T}\!\!\!\int_{Q_T} \Psi \, I^{''}_\dt \, (\ups_2(\bar{X}) -
        \ups_1(X)) \, (\ptl_y \, \ups_1(X) \, \ptl_y \, u_1(X) +
        \ptl_{\bar{y}} \, \ups_2(\bar{X}) \, \ptl_{\bar{y}} \, u_2(\bar{X}))
        \, dX\,d\bar{X} \\
        & & \quad - \dps \lim_{\dt \ra 0^+}  \int_{Q_T}\!\!\!\int_{Q_T} I^{''}_\dt \, (\ups_1(X) -
        \ups_2(\bar{X})) \, \left( \frac{\ups_2(\bar{X})}{\ups_1(X)} + \frac{\ups_1(X)}{\ups_2(\bar{X})}
        \right) \\
        & & \qquad \qquad \qquad \qquad \dps \frac{1}{2}
        \left[ (\ptl_y \, l_n \, u_1(X))^2 + (\ptl_{\bar{y}} \,
        l_n \, u_2(\bar{X}))^2 \right] \Psi \, dX\,d\bar{X} \\
        & = & \dps \lim_{\dt \ra 0^+}  \int_{Q_T}\!\!\!\int_{Q_T} I^{''}_\dt \, (\ups_1(X) -
        \ups_2(\bar{X})) \, \left(- \frac{\Psi}{2} \right)
        \frac{(\ups_1(X)-\ups_2(\bar{X}))^2}{\ups_1(X) \,
        \ups_2(\bar{X})} \\
        & & \qquad \qquad \qquad \qquad \left[ (\ptl_y \, l_n \, u_1(X))^2 + (\ptl_{\bar{y}} \,
        l_n \, u_2(\bar{X}))^2 \right] dX\,d\bar{X} \\
        & = & 0 \, ,
        \end{array}
        \end{eqnarray}
        where in the last step, we have used the definition of
        $I_\dt$ and the Lebesgue's dominant convergence theorem.

        Now, passing to the limit $h \ra 0^+$ in (3.16), and using
        the estimates (3.18), (3.26) - (3.29), (3.30) - (3.34), we
        have derived the desired estimate (3.12).  This completes
        the proof of Proposition 3.4. \qquad \qquad \qquad \qquad
        \qquad \qquad \qquad \qquad \qquad \qquad
        \qquad \qquad \qquad \qquad \qquad \qquad
        \qquad \qquad \qquad \qquad $\square$

    \end{enumerate}
\end{enumerate}

Now, we are in the position to prove the uniqueness and continuous dependence of weak
solutions to the initial-boundary value problem (2.1).

\vskip 1cm

\begin{proposition}
{\rm Let $u_1(x,y,t)$ and $u_2(x,y,t)$ be two weak solutions
to the problem (2.1) with corresponding initial data
$(u_{10}(x,y), u_{20}(x,y))$ and boundary data $(u_{11}(y,t),
u_{21}(y,t))$ and $(\ups_{01}(x,t), \ups_{02}(x,t))$ respectively. Then it holds that for almost all $t \in [0,T]$,
\begin{eqnarray}
\begin{array}{lcl}
& & \dps \int^L_0\!\!\!\int^1_0 |u_1(x,y,t) - u_2(x,y,t)| dy\,dx \\
& \leq & \dps c_6 \left\{ \int^1_0\!\!\!\int^L_0 |u_{10} (x,y) -
u_{20}(x,y)| dx\,dy + \int^t_0\!\!\!\int^1_0 |u_{11} (y,s) -
u_{21}(y,s)| dy\,ds \right. \\
& & \qquad \qquad \dps + \left. \int^t_0\!\!\!\int^L_0 |\ups_{01}
(x,s) - \ups_{02}(x,s)| dx\,ds \right\} \, ,
\end{array}
\end{eqnarray}
where $c_6 = c_6(Q,T)$ is a positive constant. }
\end{proposition}

\vskip 1cm

\begin{enumerate}
    \item[]
    \begin{enumerate}
        \item[{\bf Proof:}] This proposition follows from
        Proposition 3.4 by choosing appropriate test function
        $\phi$. Indeed, let $\rho = \rho(\xi)$ be the mollifying
        function defined in the proof of Proposition 3.4.  Fix any
        $t \in (0,T]$.  For $\var \in (0, \frac{1}{4} \min
        (t,L,1))$, we set
        $$
        \phi_{0\var} (\tau) = \int^{\tau-2\var}_{\tau-t+2\var}
        \rho_\var (s) ds, \quad \phi_{1\var} (x) = \int^{x-2\var}_{x-l+2\var}
        \rho_\var (s) ds, \quad \phi_{2\var} (y) = \int^{y-2\var}_{y-1+2\var}
        \rho_\var (s) ds \, .
        $$

        Then, clearly, $\phi_{0\var} \in C^\infty_0 (0,t)$, $\phi_{1\var} \in C^\infty_0
        (0,L)$, $\phi_{2\var} \in C^\infty_0 (0,1)$, and
        \begin{eqnarray}
        0 \leq \phi_{0\var} (\tau), \phi_{1\var} (x), \phi_{2\var}
        (y) \leq 1, \quad \forall (t,x,y) \in (0,t) \times [0,L]
        \times [0,1] \, .
        \end{eqnarray}

        Set,
        \begin{eqnarray}
        \phi_\var (x,y,\tau) = (1-y)^2 \phi_{0\var} (\tau) \, \phi_{1\var} (x)
        \, \phi_{2\var} (y), \quad \forall (x,y,\tau) \in \, Q_t \, .
        \end{eqnarray}

        Then,
        \begin{eqnarray}
        \phi_\var \, \in \, c^\infty_0 (Q_t), \quad {\rm and} \quad
        0 \leq \phi_\var (x,y,\tau) \leq 1, \quad \forall (x,y,\tau) \in \, Q_t \, .
        \end{eqnarray}

        It then follows from (3.12) in Proposition 3.4 that
        \begin{eqnarray}
        \begin{array}{lcl}
        & & I(u_1,u_2,\phi_\var) \\
        & \equiv & \dps \int_{Q_t} \{|
        \ups_1(x,y,\tau) - \ups_2(x,y,\tau)|
        (\ptl_{{t}} \, \phi_\var + \ptl_x (a\,\phi_\var) +
        \ptl_y (b\,\phi_\var) + c\,\phi_\var ) - \ptl_y \,
        \phi_\var \, \ptl_y |u_1 - u_2| \} dx\,dy\,d\tau \, \, \\
        & \geq & 0
        \end{array}
        \end{eqnarray}
        with $\ups_i=(u_i)^{-1}$.  We can rewrite
        $I(u_1,u_2,\phi_\var)$ as
        \begin{eqnarray}
        \begin{array}{lcl}
        \dps I(u_1,u_2,\phi_\var) & = & \dps \int_{Q_t} (1-y)^2 |\ups_1-\ups_2| \phi^{'}_{0\var}
        (\tau) \, \phi_{1\var} (x) \, \phi_{2\var} (y) \,
        dx\,dy\,d\tau \\
        & & \quad + \dps \int_{Q_t} (1-y)^2 a |\ups_1-\ups_2| \phi_{0\var}
        (\tau) \, \phi^{'}_{1\var} (x) \, \phi_{2\var} (y) \,
        dx\,dy\,d\tau \\
        & & \quad + \dps \int_{Q_t} b |\ups_1-\ups_2| \phi_{0\var}
        (\tau) \, \phi_{1\var} (x) \, \left( (1-y)^2 \phi_{2\var} (y) \right)_y \,
        dx\,dy\,d\tau \\
        & & \quad - \dps \int_{Q_t} \left( \ptl_y |u_1-u_2| \right) \phi_{0\var}
        (\tau) \, \phi_{1\var} (x) \, \left( (1-y)^2 \phi_{2\var} (y) \right)_y \,
        dx\,dy\,d\tau \\
        & & \quad + \dps \int_{Q_t} (\ptl_x \, a + \ptl_y \, b - c) |\ups_1-\ups_2|
        \phi_\var \, dx\,dy\,d\tau \\
        & \equiv & \dps \sum^5_{i=1} K^\var_i
        \end{array}
        \end{eqnarray}

        Each term on the right hand above can be estimated as
        follows.  First, since
        \begin{eqnarray*}
        \begin{array}{lcl}
        K^\var_1 & \equiv & \dps \int_{Q_t} (1-y)^2 |\ups_1 - \ups_2|
        \phi^{'}_{0\var} (\tau) \, \phi_{1\var} (x) \,
        \phi_{2\var} (y) \, dx\,dy\,d\tau \\
        & = & \dps \int_{Q_t} (1-y)^2 |\ups_1 - \ups_2| \left(
        \rho_\var (\tau - 2\var) - \rho_\var (\tau - t + 2\var)
        \right) \phi_{1\var} (x) \, \phi_{2y} (y) \, dx\,dy\,d\tau
        \, ,
        \end{array}
        \end{eqnarray*}
        it follows that as $\var \ra 0^+$
        \begin{eqnarray}
        K^\var_1 \lra - \left. \int^L_0\!\!\!\int^1_0 (1-y)^2
        |\ups_1 - \ups_2| \right|^{\tau=t}_{\tau=0} dy\,dx \equiv
        K_1 \, .
        \end{eqnarray}

        Similarly, as $\var \ra 0^+$,
        \begin{eqnarray}
        \begin{array}{lcl}
        K^\var_2 & \equiv & \dps \int_{Q_t} (1-y)^2 |\ups_1 - \ups_2|
        \phi_{0\var} (\tau) \, \phi^{'}_{1\var} (x) \,
        \phi_{2\var} (y) \, dx\,dy\,d\tau \\
        & & \qquad \lra \dps - \left. \int^t_0\!\!\!\int^1_0 (1-y)^2 a |\ups_1 - \ups_2|
        \right|^{x=L}_{x=0} dy\,d\tau \\
        & \equiv & K_2 \, ,
        \end{array}
        \end{eqnarray}
        and
        \begin{eqnarray}
        \begin{array}{lcl}
        K^\var_3 & \equiv & \dps \int_{Q_t} b |\ups_1 - \ups_2|
        \phi_{0\var} (\tau) \, \phi_{1\var} (x) \, \left( (1-y)^2
        \phi_{2\var} (y) \right)_y \, dy\,dx\,d\tau \\
        & = & \dps \int_{Q_t} b |\ups_1 - \ups_2| \phi_{0\var} (\tau) \, \phi_{1\var} (x) \,
        \left[ 2(y-1) \, \phi_{2\var} (y) \right. \\
        & & \dps \qquad  + \left. (1-y)^2 \left( \rho_\var(y-2\var) - \rho_\var (y-1+2\var)
        \right) \right] \, dy\,dx\,d\tau \\
        & \lra & \dps \int_{Q_t} b|\ups_1 - \ups_2| \,
        2(y-1) dx\,dy\,d\tau + \left. \int^t_0\!\!\!\int^L_0 (\g b \,
        \g|\ups_1 - \ups_2|) \right|_{y=0} dx\,dt \\
        & \equiv & K_3 \, ,
        \end{array}
        \end{eqnarray}
        here and what follows, by $\g w|_\te$, we mean the trace of
        $w$ on $\te$.

        Next, we deal with $K^\var_4$.  Note that
        \begin{eqnarray*}
        \begin{array}{lcl}
        K^\var_4 & \equiv & \dps - \int_{Q_t} (\ptl_y |u_1 - u_2|)
        \phi_{0\var} (\tau) \, \phi_{1\var} (x) \left( (1-y)^2
        \phi_{2\var} (y) \right)_y dx\,dy\,d\tau \\
        & = & \dps - \int_{Q_t} (\ptl_y |u_1 - u_2|)
        \phi_{0\var} (\tau) \, \phi_{1\var} (x) (1-y)^2 \left(
        \rho_\var (y-2\var) - \rho_\var (y-1+2\var) \right)
        dx\,dy\,d\tau \\
        & & \qquad \dps - \int_{Q_t} (\ptl_y |u_1 - u_2|)
        \phi_{0\var} (\tau) \, \phi_{1\var} (x) \, 2(y-1)
        \phi_{2\var} (y) dx\,dy\,d\tau \\
        & \equiv & K^\var_{41} + K^\var_{42} \, .
        \end{array}
        \end{eqnarray*}

        It is clear that
        $$
        K^\var_{41} \lra - \left. \int^t_0\!\!\!\int^L_0 \g
        (\ptl_y|u_1 - u_2|) \right|_{y=0} dx\,d\tau \quad {\rm as}
        \quad \var \ra 0^+ \, .
        $$

        And since,
        \begin{eqnarray*}
        \begin{array}{lcl}
        K^\var_{42} & \equiv & - \dps \int_{Q_t} (\ptl_y|u_1 -
        u_2|) \phi_{0\var} (\tau) \, \phi_{1\var} (x) \, 2(y-1) \,
        \phi_{2\var} (y) dx\,dy\,dt \\
        & = & \dps + \int_{Q_t} |u_1 - u_2| \phi_{0\var} (\tau) \,
        \phi_{1\var} (x) \big( 2\phi_{2\var} (y) + 2(y-1)
        \left(\rho_\var (y-2\var) - \rho_\var(y-1+2\var)\right)
        \big) dy\,dx\,d\tau \, ,
        \end{array}
        \end{eqnarray*}
        so as $\var \ra 0^+$,
        $$
        K^\var_{42} \lra 2 \int_{Q_t} |u_1 - u_2|
        dy\,dx\,d\tau - \left. 2 \int^t_0\!\!\!\int^L_0  \g|u_1 - u_2|
        \right|_{y=0} dx\,dt \, .
        $$

        Consequently, as $\var \ra 0^+$,
        \begin{eqnarray}
        \begin{array}{lcl}
        K^\var_4 & \lra & \dps 2\int_{Q_t} |u_1-u_2|
        dx\,dy\,d\tau - \left. 2\int^t_0\!\!\!\int^L_0 \g|u_1-u_2|
        \right|_{y=0} dx\,dt - \left. \int^t_0\!\!\!\int^L_0 \g(\ptl_y|u_1-u_2|)
        \right|_{y=0} dx\,dt \, \\
        & \equiv & K_4 \,
        \end{array}
        \end{eqnarray}

        Finally, one has that as $\var \ra 0^+$,
        \begin{eqnarray}
        \begin{array}{lcl}
        K^\var_5 & \equiv & \dps \int_{Q_t} (\ptl_x a + \ptl_y b
        +c)|\ups_1-\ups_2| \phi_\var \, dx\,dy\,d\tau \\
        & & \qquad \dps \lra
        \int_{Q_t} (\ptl_x a + \ptl_y b
        +c)|\ups_1-\ups_2| (1-y)^2 dx\,dy\,d\tau \\
        & = & K_5 \, .
        \end{array}
        \end{eqnarray}

        It follows from (3.39) - (3.45) that
        \begin{eqnarray}
        \begin{array}{lcl}
        0 & \leq & K_1+K_2+K_3+K_4+K_5 \\
        & = & \dps - \left\{ \left. \int^L_0\!\!\!\int^1_0 (1-y)^2
        \g|\ups_1 - \ups_2| \right|_{\tau=t} dx\,dy + \left. \int^t_0\!\!\!\int^1_0 (1-y)^2
        (a \g|\ups_1 - \ups_2|) \right|_{x=L} dy\,d\tau \right. \\
        & & \dps \qquad \left. + 2 \left. \int^t_0\!\!\!\int^L_0 \g|u_1 - u_2| \right|_{y=0}
        dx\,d\tau \right\} + \left. \int^L_0\!\!\!\int^1_0 (1-y)^2
        \g|\ups_1 - \ups_2| \right|_{\tau=0} dx\,dy \\
        & & \dps + \left. \int^t_0\!\!\!\int^1_0 (1-y)^2
        (a \g|\ups_1 - \ups_2|) \right|_{x=0} dy\,dt + \int_{Q_t}
        b \, 2(y-1) |\ups_1 - \ups_2| dx\,dy\,d\tau \\
        & & \dps + 2\int_{Q_t} |u_1 - u_2| dx\,dy\,d\tau +
        \int_{Q_t} (\ptl_x a + \ptl_y b + c) |\ups_1-\ups_2|
        (1-y)^2 dx\,dy\,d\tau \\
        & & \dps + \left. \int^t_0\!\!\!\int^L_0 (b \, \g|\ups_1 -
        \ups_2| - \g(\ptl_y |u_1 - u_2|)) \right|_{y=0} dx\,d\tau
        \end{array}
        \end{eqnarray}

        On the other hand, it follows from the boundary condition in
        (2.1) that
        $$
        \left( b\, \g|\ups_1-\ups_2| - \g(\ptl_y |u_1-u_2|)
        \right) \big|_{y=0} = \overline{-{\rm sgn} (u_1-u_2)} \,
        (\ups_{01} - \ups_{02}) \qquad {\rm a.e.}
        $$

        Thus,
        \begin{eqnarray}
        \left. \int^t_0\!\!\!\!\int^L_0 \left( b \,
        \g|\ups_1-\ups_2|-\g(\ptl_y|u_1-u_2|) \right) \right|_{y=0}
        dx\,dt \leq \int^t_0\!\!\!\int^L_0
        |\ups_{01}-\ups_{02}|dx\,d\tau \, .
        \end{eqnarray}

        Noting that
        \begin{eqnarray*}
        \begin{array}{ccc}
        & c_0|u_1-u_2| \leq (1-y)^2 |\ups_1-\ups_2| \leq c^{-1}_0
        |u_1-u_2| ,& {\rm and} \\
        & \big| b(x,y,t) (1-y) |\ups_1-\ups_2| \big| \leq c^{-1}_0
        |u_1-u_2| ,&
        \end{array}
        \end{eqnarray*}
        one gets from (3.46) and (3.47) that
        \begin{eqnarray}
        \begin{array}{lcl}
        & & \dps \int^L_0\!\!\!\int^1_0 |u_1(x,y,t)-u_2(x,y,t)| dy\,dx +
        \int^t_0\!\!\!\int^1_0 |u_1(L,y,t)-u_2(L,y,t)| dy\,dt \\
        & & \dps \qquad + \int^t_0\!\!\!\int^L_0 |u_1(x,0,t)-u_2(x,0,t)|
        dx\,dt \\
        & \leq & \dps \tilde{c}_0 \left\{ \int^L_0\!\!\!\int^1_0 |u_{10}(x,y)-u_{20}(x,y)|
        dy\,dx + \int^t_0\!\!\!\int^1_0 |u_1(0,y,t)-u_2(0,y,s)| dy\,dt \right. \\
        & & \dps \qquad + \left. \int^t_0\!\!\!\int^L_0
        |\ups_{01}-\ups_{02}|(x,t) dx\,dt + \int_{Q_t} |u_1-u_2|
        dx\,dy\,d\tau \right\} \\
        \end{array}
        \end{eqnarray}

        This yields the desired estimate (3.35) immediately.  So
        the proof of Proposition 3.5 is completed.
        \begin{flushright}
        $\square$
        \end{flushright}
    \end{enumerate}
\end{enumerate}

As a consequence of Proposition 3.5, we conclude that

\vskip 1cm

\begin{theorem}
{\rm The initial-boundary value problem (2.1) has a unique weak
solution $u$ as defined in Theorem 2.1.  Furthermore, such a weak
solution depends on its initial and boundary data continuously in
$L^1$-norm. }
\end{theorem}

In next section, we will study the regularity of the weak solution
to the initial-boundary value problem (2.1).

\vskip 2cm

\section{Interior Regularity}
\indent

\setcounter{equation}{0}
\setcounter{theorem}{0}

We now turn to establish the interior smoothness of the unique weak solution obtained in the previous sections. The main result of this section can be stated as follows:

\begin{theorem}
Under the same assumptions as in Theorem 1.1, the unique weak solution to (1.7) is smooth in the interior of $\Omega_T$.
\end{theorem}

The proof of the theorem is based on studies of the regularity theory for a class of ultra-parabolic equations with rough coefficients, which is rather subtle and complicated.  To illustrate our main ideas and make the presentation clear, we will only present the proof of the desired result for the case that the corresponding Euler flow is uniform, e.g.,
\begin{eqnarray}
U(x,t)\equiv 1
\end{eqnarray}

In this case, the problem (2.1) becomes
\begin{eqnarray}
\left\{ \begin{array}{ll}
\partial_t u-u^2\partial^2_y u+y\partial_x u=0 & \\
u(x,y,t=0)=u_0(x,y), & u(x=0,y,t)=u_1(y,t)\\
u(x,y=1,t=0)=0, & \partial_y u(x,y=0,t)=U_0(x,t)
\end{array}
\right.
\end{eqnarray}

It follows from Theorem 4.1 that (4.2) has a uniqueness weak solution $u\in BV(Q_T)\cap L^\infty(Q_T)$ with the properties in (2.3)-(2.8). We will study the interior regularity in $$Q_T=\{(x,y,t)| 0<t<T, 0<x<L, 0<y<1\}.$$  For any fixed point $(x_0,y_0,t_0)\in Q_T$, there exists a positive constant $\delta>0$ such that
\begin{eqnarray}
\bar{B}_\delta=\bar{B}_\delta (x_0,y_0,t_0)\triangleq \{(x,y,t)| \, |x-x_0| \leq \delta^3, |y-y_0|\leq \delta, |t-t_0|<\delta^2\} \subset Q_T
\end{eqnarray}
One can study the regularity of $u(x,y,t)$ in $B_\delta$. To this end, one can set $w(x,y,t)=u^{-1}(x,y,t)$. Then (4.2) implies that
\begin{eqnarray}
\partial_t w-\partial_y(a\partial_y w)+ y\partial_x w=0 \quad \text{on} \quad B_\delta.
\end{eqnarray}
with
\begin{eqnarray}
a=u^2(x,y,t)
\end{eqnarray}

By shifting and rescaling the independent variables
\begin{eqnarray}
(\tilde{x}, \tilde{y}, \tilde{t})=\left(\frac{x-x_0}{\delta^3}, \frac{y-y_0}{\delta}, \frac{t-t_0}{\delta^2}\right),
\end{eqnarray}
one obtains from (4.4) that
\begin{eqnarray}
\partial_{\tilde{t}} w-\partial_{\tilde{y}}(a\partial_{\tilde{y}} w) + (\tilde{y}+y_0)\partial_{\tilde{x}}w=0 \quad \text{on} \quad \tilde{B}_1(0,0,0)
\end{eqnarray}
with
\begin{eqnarray}
\tilde{B}_1(0,0,0)=\{(\tilde{x}, \tilde{y}, \tilde{t})| \, |\tilde{x}|<1, |\tilde{y}|<1, |\tilde{t}|<1\}
\end{eqnarray}

The equation (4.7) can be simplified further by introduce the following transformation
\begin{eqnarray}
\tau=\tilde{t}, \quad \xi=\tilde{x}-y_0\,\tilde{t}, \quad \eta=\tilde{y},
\end{eqnarray}
such that
\begin{eqnarray}
\partial_\tau w-\partial_\eta (a\partial_\eta w)+\eta\partial_\xi w=0 \quad \text{on} \quad \tilde{\tilde{B}}_1 (0,0,0)
\end{eqnarray}
with $\tilde{\tilde{B}}_1 (0,0,0)=\{(\tau,\xi,\eta)| \, |\xi+y_0\tau|<1, |\eta|<1, |\tau|<1\}$.

Choose constant $\delta_1>0$ so that
\begin{eqnarray}
\bar{B}_{\delta_1} (0,0,0)\equiv\{ (\tau,\xi,\eta)| \, |\xi|\leq \delta^3_1, |\eta|\leq\delta_1, |\tau|^2\leq \delta^2_1\}\subset \tilde{\tilde{B}}_1 (0,0,0)
\end{eqnarray}

Rescale again by
\begin{eqnarray}
(\tilde{\xi},\tilde{\eta},\tilde{\tau})=\left( \frac{\xi}{\delta^3_1}, \frac{\eta}{\delta_1}, \frac{\tau}{\delta^2_1}\right)
\end{eqnarray}
so that one can get from (4.10) that
\begin{eqnarray}
\partial_{\tilde{\tau}} w-\partial_{\tilde{\eta}} (a\partial_{\tilde{\eta}}w)+ \tilde{\eta}\partial_{\tilde{\xi}}w=0 \quad \text{on} \quad \tilde{B}_1 (0,0,0)
\end{eqnarray}
with
\begin{eqnarray}
\tilde{B}_1(0,0,0)=\{(\tilde{\xi},\tilde{\eta},\tilde{\tau})| \, |\tilde{\xi}|<1, |\tilde{\eta}|<1, |\tilde{\tau}|<1 \}
\end{eqnarray}

Thus we will study the equation (4.13). For notational convenience, we will use $u$ for $w$, $(x,y,t)$ for $(\tilde{\xi},\tilde{\eta},\tilde{\tau})$, and $B_1$ for $\tilde{B}_1$. Hence consider
\begin{eqnarray}
\partial_t u-\partial_y(a\,\partial_y u)+y\partial_x u=0 \quad \text{on} \quad B_1(0,0,0)
\end{eqnarray}
with the assumption that
\begin{eqnarray}
a\in L^\infty, \Lambda^{-1} \leq a\leq \Lambda \quad \text{on} \quad B_1(0,0,0)
\end{eqnarray}
for a positive constant $\Lambda$.

We will modify the idea of Krushkov in [5] on the $C^\alpha$-regularity theory of weak solutions to uniform parabolic equations to analyze the regularity of weak solutions to the ultra parabolic equation (4.15) with the assumption (4.15). In particular, we introduce a "weak" form of Poincar\'e inequality. The key step is to establish suitable oscillation estimates. To this end, one needs the following notations. For any given constants $r\in(0,1], \alpha, \beta\in(0,1)$, set
\begin{eqnarray*}
& B_1=B_r(0,0,0)=\{(x,y,t)| \, |x|<r^3, |y|<r, |t|<r^2\},\\
& B^\pm_r\triangleq B_r\cap\{t\gtrless 0\},\\
& C_r=\{(x,y)| \, |x|<r^3, \, |y|<r\},\\
& H_{t,h}\triangleq H^u_{t,h}=\{(x,y)| (x,y)\in C_{\beta r}, u(x,y,t)\geq h\}
\end{eqnarray*}
for any given weak solution $u$ to (4.15), and any fixed $t\in(-r^2,0)$ and $h\in(0,1)$. Then the following estimates holds.

\begin{proposition}
Let $u(x,y,t)$ be a non-negative weak solution to (4.15) in $B^-_1$ with the property that
\begin{eqnarray}
\text{mes}\,\{(x,y,t)\in B^-_r, u(x,y,t)\geq 1\}\geq \frac{1}{2}\,\text{mes}\, B^-_r.
\end{eqnarray}
\end{proposition}

Then there exist constants $\alpha, \beta, h_1\in(0,1)$, $\alpha\ll 1$, $\beta\approx 1$, depending only on $\Lambda$, such that for almost all $t\in(-\alpha r^2,0)$, $0<h\leq h_1$, it holds that
\begin{eqnarray}
\text{mes}\, H_{t,h}\geq\frac{1}{11}\,\text{mes}\, C_{\beta r}
\end{eqnarray}

{\bf Proof}: For $h\in(0,\frac{1}{2})$, set
\begin{eqnarray}
V(x,y,t)=ln^+\frac{1}{u(x,y,t)+h^{\frac{9}{8}}}=(-ln(u+h^{\frac{9}{8}})^+,
\end{eqnarray}
where, and in the following, $f^+$ denote the positive part of the function $f$.

Then it can be checked easily that $V$ is a non-negative weak solution to the following equation
\begin{eqnarray}
\partial_t V-\partial_y(a\partial_y V)+a(\partial_y V)^2+y\partial_x V=0\quad \text{on}\quad B^-_1.
\end{eqnarray}

Let $\chi(s)$ be a smooth cut-off function such that
\begin{eqnarray}
\left\{ \begin{array}{ll}
\chi\in C^\infty([0,\infty)), 0\leq \chi \leq 1 \quad \text{for}\quad 0\leq s\leq \beta r, \chi(s)=0\quad\text{for}\quad s\geq r & \\
|\chi'(s)|\leq \frac{2}{(1-\beta)r} \quad \text{for all}\quad s\geq0 &
\end{array}
\right.
\end{eqnarray}

It then follows from (4.20) that for almost all $t,\tau \in(0-r^2,0)$, $\tau\leq t$, it holds that
\begin{eqnarray}
\begin{array}{cl}
& \displaystyle \int_{C_r} \chi^2(|y|) V(x,y,t)dxdy + \int^t_\tau\!\!\int_{C_r} |\partial_y V|^2 a\chi^2(|y|)dxdyds\\
= & \displaystyle \int_{C_r} \chi^2 (|y|)V(x,y,\tau)dxdy - \int^t_\tau\!\!\int_{C_r} (2\chi\partial_y\chi)a\partial_y V dxdyds - \int^t_\tau\!\!\int^r_{-r}\chi^2 yV|^{r^3}_{-r^3}dyds
\end{array}
\end{eqnarray}

Cauchy-Schwartz inequality yields
\begin{eqnarray*}
\begin{array}{lcl}
\displaystyle -\int^t_\tau\!\!\int_{C_r}(2\chi\partial_y \chi)a\partial_y V dxdyds & \leq & \displaystyle \frac{1}{2}\int^t_\tau\!\!\int_{C_r} a\chi^2 |\partial_y V|^2dxdyds + 4\int^t_\tau\!\!\int_{C_r}a|\partial_y\chi|^2 dxdyds\\
& \leq & \displaystyle \frac{1}{2}\int^t_\tau\!\!\int_{C_r} a\chi^2 |\partial_y V|^2dxdyds + \frac{16\Lambda}{(1-\beta)^2 r^2}(t-\tau)\text{mes}\, C_r\\
& \leq & \displaystyle \frac{1}{2}\int^t_\tau\!\!\int_{C_r} a\chi^2 |\partial_y V|^2 dxdyds + \frac{16\Lambda}{\beta^4(1-\beta)^2}\text{mes}\, C_{\beta r}.
\end{array}
\end{eqnarray*}

Due to the definition of $V$, one gets by direct calculations that
$$-\int^t_\tau\!\!\int^r_{-r}\chi^2 yV\big|^{r^3}_{-r^3} dyds \leq (t-\tau) r^2 lnh^{-\frac{9}{8}}\leq \frac{1}{4}\beta^{-4}lnh^{-\frac{9}{8}}\,\text{mes}\, C_{\beta r}$$
These together with (4.22) yield
\begin{eqnarray}
\begin{array}{ll}
& \displaystyle \int_{C_r} \chi^2(|y|) V(x,y,t) dxdy + \frac{1}{2}\int^t_\tau\!\!\int_{C_r} a\chi^2|\partial_y V|^2 dxdyds\\
\leq & \displaystyle \int_{C_r} \chi^2(|y|)V(x,y,\tau)dxdy + \frac{1}{4}\beta^{-4} lnh^{-\frac{9}{8}} \text{mes} \, C_{\beta\gamma}+\frac{16\Lambda}{\beta^4(1-\beta)^2}\text{mes} \, C_{\beta\gamma}.
\end{array}
\end{eqnarray}

To estimate the last term on the right hand side above, one can define
\begin{eqnarray}
\mu(t)=\text{mes}\, \{(x,y)\in C_r, \, u(x,y,t) \geq 1\} \quad \text{for any} \quad t\in[-r^2,0].
\end{eqnarray}

It then follows from the assumption (4.17) and the definition (4.24) that
\begin{eqnarray*}
\int^{\circ}_{-r^2}\mu(t)d\tau=\text{mes}\, \{(x,y,t)\in B^-_r, \, u(x,y,t)\geq 1\} \geq \frac{1}{2}\, \text{mes} \, B^-_r=\frac{1}{2}r^2\, \text{mes}\, C_r
\end{eqnarray*}
which implies that for any constant $\alpha\in (0,\frac{1}{2})$,
$$\int^{-\alpha r^2}_{-r^2} \mu(t) dt\geq \frac{1}{2} r^2\, \text{mes}\, C_r-\int^{\circ}_{-\alpha r^2}\mu(t) dt\geq (\frac{1}{2}-\alpha)r^2\, \text{mes}\, C_r.$$

Thus there exists $\tau\in(-r^2,-\alpha r^2)$ such that
\begin{eqnarray}
\mu(\tau)\geq (\frac{1}{2}-\alpha)(1-\alpha)^{-1}\, \text{mes}\, C_r.
\end{eqnarray}

Note that the definition of $V$ and (4.25) imply that
\begin{eqnarray}
\begin{array}{rcl}
\displaystyle \int_{C_r} V(x,y,\tau)dxdy & = & \displaystyle \int_{C_r\cap\{u\leq 1\}} V(x,y,\tau)dxdy\\
& \leq & \displaystyle (lnh^{-\frac{9}{8}})\, \text{mes}\, \left\{C_r\cap\{u\leq 1\}\right\}\\
& = & \displaystyle lnh^{-\frac{9}{8}} (\text{mes}\, C_r-\mu(\tau))\leq \frac{1}{2}(1-\alpha)^{-1}\, \text{mes}\, C_r\, lnh^{-\frac{9}{8}}.
\end{array}
\end{eqnarray}

We now turn to estimate the first term on the left hand side of (4.23). Note  that for $h\in(0,\frac{1}{2})$, $V(x,y,t)\geq ln\frac{1}{h+h^{\frac{9}{8}}}$ for all $(x,y)\not\in H_{t,h}$. It holds that
\begin{eqnarray}
\int_{C_r}\chi^2 V(x,y,t)dxdy\geq \int_{C_{\beta r}}V(x,y,t)dxdy\geq ln\frac{1}{h+h^{\frac{9}{8}}}\,\text{mes}\, (C_{\beta r}\setminus H_{t,h}).
\end{eqnarray}

It follows from (4.23), (4.26) and (4.27) that
\begin{eqnarray}
\text{mes}\, (C_{\beta r}\setminus H_{t,h})\leq \frac{lnh^{-\frac{9}{8}}}{ln\frac{1}{h+h^{\frac{9}{8}}}} \left( \frac{1}{2}(1-\alpha)^{-1}\beta^{-4}+\frac{1}{4}\beta^{-4}+\frac{16\Lambda}{\beta^4(1-\beta)^2}\frac{1}{lnh^{-\frac{9}{8}}}\right)\,\text{mes}\,C_{\beta r}
\end{eqnarray}

Since $\lim_{h\rightarrow 0^+} \frac{ln h^{-\frac{9}{8}}}{ln\frac{1}{h+h^{\frac{9}{8}}}}=\frac{9}{8}$, so there exist constants $\alpha$, $\beta$, $h_1 \in(0,1)$ with property that $\alpha\ll 1$, $\beta\approx 1$, and $h_1=h_1(s)$ suitably small such that for $h\in(0,h_1)$,
$$\frac{ln\,h^{-\frac{9}{8}}}{ln\frac{1}{h+h^{\frac{9}{8}}}} \left( \frac{1}{2}(1-\alpha)^{-1} \beta^{-4}+\frac{1}{4}\beta^{-4} + \frac{16\Lambda}{\beta^4(1-\beta)^2}\frac{1}{ln\,h^{-\frac{9}{8}}}\right)\leq \frac{10}{11}.$$

This and (4.28) yield the desired estimate (4.18). The proof of Proposition 4.2 is complete.

The next key element of our analysis is a "weak" form of Poincar\'e's inequality based on the fundamental solution of the equation (4.15) with $a\equiv 1$. This kind of "weak" form of Poincar\'e's inequality is needed for non-negative sub-solutions to (4.15) defined as follows:

We denote $$L^2_+(\Omega_{x,t}; H^{-1}_y)=\{u\in BV(\Omega)| \text{for any}\quad \varphi\in C^\infty_0(\Omega), \varphi\geq 0\quad
\frac {\int_\Omega (\partial_t u+y\partial_x u)\varphi dx\,dy\,dt}{|\partial_y\varphi|_{L^2(\Omega)}+|\varphi|_{L^2(\Omega)}}< \infty\}.$$

\begin{remark} $u$ is said to be a weak subsolution to (4.15) on $\Omega$ if $\partial_t u+y\partial_x u\in L^2_+(\Omega_{x,t}; H^{-1}_y)$, $\partial_y u \in L^2_{\text{loc}}(\Omega)$, $u\in L^\infty_{\text{loc}}(\Omega)$, and for any $\varphi\in C^\infty_0(\Omega)$ with $\varphi\geq 0$, it holds that
\begin{eqnarray}
-\int_\Omega (a\partial_y \varphi\partial_y u+ (\partial_t u+y\partial_x u)\varphi) dx\,dy\,dt\geq 0.
\end{eqnarray}
\end{remark}

Consider the following basic ultra-parabolic equation
\begin{eqnarray}
\mathcal{L}_0 u\equiv \partial_t u-\partial^2_y u + y\partial_x u=0
\end{eqnarray}

Set $z=(x,y,t)$ and $\zeta=(\xi,\eta,\tau)$. Then the fundamental solution of $\mathcal{L}_0$ can be constructed in [14, 15, 10] as
\begin{eqnarray}
\Gamma_0 (z,\zeta)=\Gamma_0(\zeta^{-1}oz,0)=\left\{
\begin{array}{cc}
\frac{\sqrt{3}}{2\pi(t-\tau)^2}\exp \left\{-\frac{(y-\eta)^2}{3(t-\tau)}-\frac{3}{(t-\tau)^3} (x-\xi-\frac{(t-\tau)}{2} (y+\eta))^2\right\} & t>\tau\\
0 & t\leq\tau,
\end{array}
\right.
\end{eqnarray}
where $\zeta oz$ is the left translation of $z$ and $\zeta$ in the group law associated with $\mathcal{L}_0$. It is known that
\begin{eqnarray}
\int_{\mathbb{R}^2}\Gamma_0(z,\zeta)dx\,dy=\int_{\mathbb{R}^2} \Gamma_0(z,\zeta)d\xi\,d\eta=1 \quad \text{for} \quad t>\tau,
\end{eqnarray}
and
\begin{eqnarray}
\Gamma_0(\delta_\mu oz,0)=\mu^{-4} \Gamma_0(z,0), \quad \forall \, z\not= 0, \quad \mu>0
\end{eqnarray}
where $\delta_\mu$ is the dilations associated with $\mathcal{L}_0$. To establish a ``weak" form of the Poincar\'e's inequality for non-negative weak subsolutions to (4.15), one need to construct suitable test functions. To this end, one chooses a cut-off function $\chi\in C^\infty([0,\infty))$ with the properties that $0<\chi(s)<1$ for $0\leq s<r$, $\chi(s)=0$ for all $s>r$, $\chi(s)\equiv 1$ for $0\leq s\leq \theta^{\frac{1}{6}}\,r$, and
\begin{eqnarray}
\left\{ \begin{array}{ll}
0\leq-\chi'(s)\leq \frac{2}{(1-\theta^{\frac{1}{6}})r}, |\chi''(s)|\leq \frac{C}{r^2}, & \text{for all}\quad s\in\mathbb{R}^1_+,\\
|\chi'(s)|\geq C(\beta_1,\beta_2)r^{-1}>0 & \text{for all}\quad s\in[\beta_1r, \beta_2r], \quad \theta^{\frac{1}{6}}<\beta_1<\beta_2<1,
\end{array}
\right.
\end{eqnarray}
where $\theta\in (0,2^{-6})$ is a positive constant to be chosen and $C$ and $C(\beta_1,\beta_2)$ are fixed positive constants. Set
\begin{eqnarray}
Q^-_\theta=\left\{(x,y,t)| -r^2\leq t\leq 0, \quad |y|\leq \frac{r}{\theta}, \quad |x|\leq \frac{r^3}{\theta}\right\}
\end{eqnarray}
and
\begin{eqnarray}
\phi(x,y,t)=\phi_0(x,y,\tau)\phi_1(x,y,t)
\end{eqnarray}
with
\begin{eqnarray}
\phi_0(x,y,t)=\chi([\theta^2 x^2-6tr^4]^{\frac{1}{6}}), \phi_1(x,y,t)=\chi(\theta|y|)
\end{eqnarray}

Then the following elementary facts can be verified by direct computations:
\begin{lemma}
It holds that
\begin{enumerate}
\item
\begin{eqnarray}
(-y\partial_x-\partial_t)\phi_0(z)\leq 0 \quad \text{for}\quad z\in Q^-_\theta;
\end{eqnarray}
\item \begin{eqnarray}
\phi(z)\equiv 1 \quad \text{on}\quad B^-_{\theta r};
\end{eqnarray}
\item \begin{eqnarray}
\text{supp}\,\phi\cap\{(x,y,t)|t\leq 0\}\subset Q^-_\theta;
\end{eqnarray}
\item There exists a positive constant $\alpha_1 \in(0,\min(\alpha,\frac{1}{12}))$ such that
\begin{eqnarray}
\{(x,y,t)|-\alpha_1\,r^2\leq t\leq 0, \, (x,y)\in C_{\beta r}\}\subseteq \text{supp}\,\phi,
\end{eqnarray}
\item Assume that $\alpha_1>\theta$. Then
\begin{eqnarray}
0<\phi_0(z)<1 \quad \text{on}\quad \{(x,y,t)|-\alpha_1\,r \leq t\leq -\theta r^2, \quad (x,y)\in C_{\beta r}\}.
\end{eqnarray}
\end{enumerate}
\end{lemma}

Now we choose $\theta$ and $\alpha_1$ such that all the requirements in Lemma 4.4 are satisfied. Thus (4.38)-(4.42) hold true. Then we have the following key inequality of Poincar\'e type for any non-negative weak subsolutions to (4.15).

\begin{proposition}
Let $w$ be a non-negative weak subsolution to (4.15) in $B^-_1$. Then there exists a uniform positive constant $C$ such that for any $r\in(0,\theta)$, it holds that
\begin{eqnarray}
\int_{B^-_{r\theta}} |(w(z)-I_0)^+|^2 dz\leq C\theta^2 r^2\int_{B^-_{r/\theta}}|\partial_y w(z)|^2 dz,
\end{eqnarray}
where $\displaystyle I_0=\sup_{B^-_{\theta r}} I_1(z)$ with
\begin{eqnarray}
I_1(z)=\int_{B^-_{r/\theta}} \partial_\eta\phi(\zeta)\partial_\eta \Gamma_0(z,\zeta)w(\zeta)d\zeta+\int_{B^-_{r/\theta}}\Gamma_0(z,\zeta)|\partial_\tau\phi(\zeta)+\eta\partial_\xi\phi(\zeta))w(\zeta)d\zeta.
\end{eqnarray}
\end{proposition}

{\bf Proof}: Since $\Gamma_0(z,\zeta)$ is the fundamental solution to (4.40), so a standard approximation [10] implies that
\begin{eqnarray}
(w\phi)(z)=\int_{\mathbb{R}^{2+1}} [\partial_\eta\Gamma_0\partial_\eta(w\phi)+\Gamma_0(\partial_\tau+\eta\partial_\xi)(w\phi)]d\zeta.
\end{eqnarray}

It follows from this, Lemma 4.4 and $Q^-_\theta \subseteq B^-_{r/\theta}$ that for $z\in B^-_{\theta r}$,
\begin{eqnarray}
\begin{array}{rcl}
w(z) & = & \displaystyle \int_{B^-_{r/\theta}} [\partial_\eta \Gamma_0\partial_\eta(w\phi)+\Gamma_0(\partial_\tau+\eta\partial_\xi)(w\phi)]d\zeta\\
& = & \displaystyle I_1(z)+\int_{B^-_{r/\theta}} [-\Gamma_0\partial_\eta\phi\partial_\eta w+\phi\partial_\eta\Gamma_0\partial_\eta w+\phi\Gamma_0(\partial_\tau w+\eta\partial_\xi w)]d\zeta\\
& = & I_1(z)+I_2(z)+I_3(z),
\end{array}
\end{eqnarray}
where $I_1(z)$ is given by (4.44), while $I_2$ and $I_3$ are given by
\begin{eqnarray}
& \displaystyle I_2(z)=\int_{B^-_{r/\theta}}[-(1+a(\zeta))\Gamma_0(z,\zeta)\partial_\eta\phi(\zeta)\partial_\eta w(\zeta)+(1-a(\zeta))\phi(\zeta)\partial_\eta\Gamma_0(z,\zeta)\partial_\eta w(\zeta)]d\zeta,\\
& \displaystyle I_3(z)=\int_{B^-_{r/\theta}} [a(\zeta)\partial_\eta w(\zeta)\partial_\eta (\Gamma_0(z,\zeta)\phi(\zeta))+\phi(\zeta)\Gamma_0(z,\zeta)(\partial_\tau w+\eta\partial_\xi w)]d\zeta.
\end{eqnarray}

First, note that by approximation if necessary, one can take $\phi(\zeta)\Gamma_0(z,\zeta)$ as a non-negative test function (4.29), thus it follows from Definition 4.3 that $I_3(z)\leq 0$ (as in the proof of Lemma2.5 in [10]). This and (4.46) imply that for any $z\in B^-_{\theta r}$,
\begin{eqnarray}
0\leq (w(z)-I_0)^+\leq (w(z)-I_1(z))^+\leq I_2(z)^+.
\end{eqnarray}

It remains to estimate $I_2(z)\equiv I_{21}(z)+I_{22}(z)$ as follows.
Note that for $z\in B^-_{\theta r}$,
\begin{eqnarray*}
\begin{array}{rcl}
I_{22}(z) & = & \displaystyle \int_{B^-_{r/\theta}}((1-a(\zeta)\phi(\zeta)\partial_\eta w(\zeta))\partial_\eta\Gamma_0(z,\zeta)d\zeta\\
& = & \displaystyle \int_{\mathbb{R}^3} [(1-a(\zeta))\phi(\zeta)\chi(\tau)]\partial_\eta w(\zeta)\partial_\eta\Gamma_0(z,\zeta)d\zeta
\end{array}
\end{eqnarray*}
where $\chi_{\tau\leq 0}$ is the characteristic function on the set $\{\tau\leq 0\}$. It follows from (4.40) and Pascucci-Polidoro's estimate (Corollary 2.2 in [10]) that
\begin{eqnarray*}
\begin{array}{rcl}
\displaystyle ||I_{22}||_{L^3(B^-_{\theta r})} & \leq & C||[(1-a(\zeta))\chi_{\tau\leq 0}\phi(\zeta)]\partial_\eta w(\zeta)||_{L^2(\mathbb{R}^3)}\\
& \leq & C||\partial_\eta w||_{L^2(B^-_{r/\theta})},
\end{array}
\end{eqnarray*}
which implies that
\begin{eqnarray}
||I_{22}||_{L^2(B^-_{\theta r})}\leq C\theta r||\partial_\eta w||_{L^2(B^-_{r/\theta})}.
\end{eqnarray}

Next, note that for $z\in B^-_{\theta r}$,
\begin{eqnarray*}
\begin{array}{rcl}
I_{21}(z) & = & \displaystyle -\int_{B^-_{r/\theta}}(1+a)\Gamma_0(z,\zeta)\partial_\eta w(\zeta)\partial_\eta\phi(\zeta)d\zeta\\
& = & \displaystyle -\int_{\mathbb{R}^{2+1}} \Gamma_0(z,\zeta)((1+a)\partial_\eta w(\zeta)\partial_\eta\phi(\zeta)\chi_{\tau\leq 0})d\zeta.
\end{array}
\end{eqnarray*}

It then follows from the Pascucci-Polidoro estimate [10], (4.40) and (4.34) that
\begin{eqnarray*}
\begin{array}{rcl}
||I_{21}||_{L^6(B^-_{\theta r})} & \leq & ||\Gamma_0((1+a)\partial_\eta w\partial\eta\phi\chi_{\tau\leq 0})||_{L^6(\mathbb{R}^{2+1})}\\
& \leq & C||(1+a)\chi_{\tau\leq 0}\partial_\eta w\partial_\eta\phi||_{L^2(\mathbb{R}^{2+1})}\\
& \leq & C||\chi_{\tau\leq 0}\partial_\eta\phi\partial_\eta w||_{L^2(B^-_{r/\theta})}\leq \frac{C}{r}||\partial_\eta w||_{L^2(B^-_{r/\theta})},
\end{array}
\end{eqnarray*}
which yields that
\begin{eqnarray}
||I_{21}||_{L^2(B^-_{\theta r})} \leq C\theta^2 r||\partial_\eta w||_{L^2(B^-_{r/\theta})}.
\end{eqnarray}

Thus, the desired estimate (4.43) follows from (4.49), (4.50) and (4.51). This completes the proof of Proposition 4.5.

Proposition 4.2 will be applied to a special class of non-negative weak solutions of (4.15). Indeed, let $u$ be a non-negative weak solution to (4.15) in $B^-_1$. Set
\begin{eqnarray}
w(z)=ln^+ \frac{h}{h^{\frac{9}{8}}+u(z)}.
\end{eqnarray}

Then it can be checked that $w(z)$ is a non-negative weak subsolution to (4.15). We will apply Proposition 4.5 to estimate $w$. To this end, one needs to derive a key estimate on $I_0$ which is a ``mean value" of $w$.

\begin{lemma}
Let $u$ be a non-negative weak solution to (4.15) in $B^-_1$ satisfying the assumptions in Proposition 4.1, and $w(z)$ be defined by (4.52). Then there exist uniform positive constants $h_0(\leq h_1)$, $\theta$, $\lambda_0$ (independent of $u$) such that
\begin{eqnarray}
\lambda_0<1, \quad |I_0|\leq \lambda_0\, ln(h^{-\frac{1}{8}}) \quad \text{for all} \quad r<\theta, \, 0<h\leq h_0,
\end{eqnarray}
with $I_0$ defined in Proposition 4.5.
\end{lemma}

{\bf Proof}: Since $u$ is a non-negative weak solution to (4.15), it can be checked directly that $w(z)$ given by (4.52) is a non-negative weak subsolution to (4.15) by Definition 4.1. Furthermore, following the same arguments in Lemma 4.4 and Proposition 4.5, one can check easily that Proposition 4.5 applies to $w$. So let $I_0$ be the corresponding ``mean value" of $w$ defined in Proposition 4.5. Then
\begin{eqnarray}
\begin{array}{rcl}
I_0 & = & \displaystyle \sup_{B^-_{\theta r}} I_1(z) \quad \text{with}\\
I_1(z) & = & \displaystyle \int_{B^-_{r/\theta}} [-\Gamma_0(z,\zeta)\partial^2_\eta \phi (\zeta)w(\zeta)]d\zeta + \int_{B^-_{r/\theta}} \Gamma_0(z,\zeta)(\partial_\tau\phi +\eta\partial_\zeta\phi)w(\zeta)d\zeta \equiv I_{11}+I_{12}
\end{array}
\end{eqnarray}
where one has integrated by parts and used Lemma 4.4 and the structure of $\partial_\eta\phi(\xi)$.

We start to estimate $I_{11}=\int_{B^-_{r/\theta}} [-\Gamma_0(z,\xi)\partial^2_\eta\phi(\xi)w(\xi)]d\xi$. Note that for any $z\in B^-_{\theta r}$, $\text{supp} ((\partial^2_\eta\phi(\zeta))\Gamma_0(z,\zeta)w(\zeta))\subset B^-_{r/\theta}$. Thus
\begin{eqnarray}
\begin{array}{rcl}
|I_{11}| & \leq & \displaystyle \int_{B^-_{r/\theta}} |\partial^2_\eta\phi(\zeta)|w(\zeta)\Gamma_0 (z,\zeta)d\zeta\\
& \leq & \displaystyle ln\, h^{-\frac{1}{8}} \int^0_{-r^2}\!\!\int_{\mathbb{R}^2} |\partial^2_\eta\phi(\zeta)|\Gamma_0(z,\zeta)d\zeta\\
& \leq & \displaystyle ln\,h^{-\frac{1}{8}} r^2 \, \sup_{\text{supp}\phi} |\partial^2_\eta\phi(\zeta)|,
\end{array}
\end{eqnarray}
where one has used (4.32) and the fact
$$B^-_{r/\theta}\cap \, \text{supp}\phi\subset Q^-_\theta$$
due to Lemma 4.4. On the other hand, it follows from (4.36)-(4.37), (4.34), and direct computatons that
$$ |\partial^2_\eta\phi(\zeta)|\leq c\theta^2 r^{-2}. $$

This, together with (4.55), yields
\begin{eqnarray}
|I_{11}|\leq c\theta^2 ln (h^{-\frac{1}{8}}).
\end{eqnarray}

We now turn to the estimate $I_{12}=\int_{B^-_{r/\theta}} \Gamma_0(z,\zeta)(\partial_\tau\phi(\zeta)+\eta\partial_\xi\phi(\zeta))w(\zeta)d\zeta$.

First, note that for suitably small $\theta>0$, $z\in B^-_{r/\theta}$, $\phi(z)=1$. Thus for $z\in B^-_{r/\theta}$,
\begin{eqnarray*}
\begin{array}{rcl}
1 & = & \displaystyle \int_{B^-_{r/\theta}} \Gamma_0(z,\zeta)(\partial_\tau+\eta\partial_\xi-\partial^2_\eta)\phi(\zeta)d\zeta\\
& = & \displaystyle \int_{B^-_{r/\theta}}\phi_1\Gamma_0(z,\zeta)(\partial_\tau\phi_0+\eta\partial_\xi\phi_0)d\zeta+\int_{B^-_{r/\theta}}\Gamma_0(z,\zeta)(-\partial^2_\eta\phi(\zeta))d\zeta.
\end{array}
\end{eqnarray*}

This and the arguments for (4.55)-(4.56) show that
\begin{eqnarray}
\int_{B^-_{r/\theta}}\phi_1(\zeta)\Gamma_0(z,\zeta)(\partial_\tau\phi_0(\zeta)+\eta\partial_\xi\phi_0(\zeta))d\zeta=1+O(1)\theta^2.
\end{eqnarray}

Next, let $\alpha_1$, $\beta$ and $h_1$ be given in Proposition 4.2 and set
\begin{eqnarray}
S=\left\{\zeta=(\xi,\eta,\tau)|-\alpha_1 r^2\leq \tau\leq -\frac{\alpha_1}{2} r^2, \quad (\xi,\eta)\in C_{\beta r}, \quad w(\zeta)=0\right\}
\end{eqnarray}

Note that if $\zeta$ is such that $u(\zeta)\geq h$, then $w(\zeta)=0$, and $u$ is assumed to satisfy the assumptions in Proposition 4.2. Thus it follows from Proposition 4.2 and the definition of $S$ that there exists a positive constant $c(\alpha_1,\beta)>0$ such that
\begin{eqnarray}
\text{mes}\, S\geq c(\alpha_1,\beta)r^6 \quad \text{for}\quad h\leq h_1
\end{eqnarray}

Furthermore, it follows from the construction of $\chi$ ((4.34)) and Lemma 4.4 that for $\theta$ suitably small and any $\zeta\in S$, $\phi_1(\zeta)=1$, $\phi_0(\zeta)>0$, $\theta^{\frac{1}{6}} r<(\theta^2 |\xi|^2-6\tau\,r^4)^{\frac{1}{6}}<r$, $6r^4-2\eta\xi\theta^2>3r^4$, and
$$|\chi'([\theta^2 |\xi|^2-6\tau r^4]^{\frac{1}{6}})|\geq \frac{c(\alpha_1,\theta)}{r}>0$$
with a uniform constant $c(\alpha_1,\theta)$ independent of $r$. Hence
\begin{eqnarray}
\begin{array}{cl}
& \displaystyle \int_S \Gamma_0(z,\zeta)(\partial_\tau\phi_0(\zeta)+\eta\partial_\xi\phi_0(\zeta))\phi_1(\zeta)d\zeta\\
= & \displaystyle \int_S \Gamma_0(z,\zeta) |\chi' ([\theta^2|\xi|^2 - 6\tau r^4]^\frac{1}{6}) \, \frac{1}{6}[\theta^2\xi^2-6\tau r^4]^{\frac{5}{6}} [6r^4-2\eta\xi\theta^2]d\zeta\\
\geq & \displaystyle c(\alpha_1,\theta) \int_S r^{-2} \Gamma_0(z,\zeta)\alpha\zeta\geq \underline{c}>0,
\end{array}
\end{eqnarray}
where $\underline{c}=c(\alpha_1,\beta,\theta)>0$, and one has used (4.59) and the fact that
$$\Gamma_0(z,\zeta)\geq cr^{-4} \quad \text{for} \quad \tau\leq-\frac{\alpha_1}{2}r^2, \quad z\in B^-_{\theta r}$$
if $\theta^2\leq \frac{1}{4}\alpha_1$. We are now ready to estimate $I_{12}$. It follows from (4.58) and Lemma 4.4 that
\begin{eqnarray}
\begin{array}{rcl}
\displaystyle |I_{12}| & = & \displaystyle \left| \int_{B^-_{r/\theta}} \Gamma_0(z,\zeta)(\partial_\tau\phi(\zeta)+\eta\partial_\xi\phi(\zeta))w(\zeta)d\zeta\right|\ \\
& = & \displaystyle \int_{B^-_{r/\theta}\setminus S}\Gamma_0(z,\zeta)(\partial_\tau\phi(\zeta)+\eta\partial_\xi\phi(\zeta))w(\zeta)d\zeta\\
& \leq & \displaystyle ln\,h^{-\frac{1}{8}} \int_{B^-_{r/\theta}\setminus S} \Gamma_0(z,\zeta)(\partial_\tau\phi(\zeta)+\eta\partial_\xi\phi(\zeta))d\zeta\\
& = & \displaystyle ln\, h^{-\frac{1}{8}} \left( \int_{B^-_{r/\theta}} \Gamma_0(z,\zeta)(\partial_\tau\phi(\zeta)+\eta\partial_\xi\phi(\zeta))d\zeta - \int_S \Gamma_0(z,\zeta)(\partial_\tau\phi(\zeta)+\eta\partial_\xi\phi(\zeta))d\zeta\right)\\
& = & \displaystyle ln\, h^{-\frac{1}{8}} \left(1+O(1)\theta^2 - \int_S \Gamma_0 (z,\zeta)(\partial_\tau\phi(\zeta)+\eta\partial_3\phi(\zeta))d\zeta\right)\\
& \leq & (1-\underline{c}+O(1)\theta^2)ln\,h^{-\frac{1}{8}},\\
\end{array}
\end{eqnarray}
where one has used (4.57) and (4.60).

Consequently, it holds that
$$ I_0=\sup_{B^-_{\theta r}} I_1(z)\leq (1-\underline{c}+O(1)\theta^2)ln\,h^{-\frac{1}{8}} $$
for any $r<\theta$, $h\leq h_1$, and $\theta$ being suitably small. Set $\lambda_0=1-\underline{c}+O(1)\theta^2$. Then the conclusion in Lemma 4.2 holds. This completes the proof of Lemma 4.6. \hskip 8cm $\square$

We are now ready to give the main step in the oscillation estimates.
\begin{lemma}
Let $u$ be a non-negative weak solution to (4.15) in $B^-_1$ satisfying the assumption in Proposition 4.2. Then there exist positive constants $\bar{h}$ and $\bar{\theta}$ in (0,1), which depend only on $\lambda_0$ and $\Lambda$, such that
\begin{eqnarray}
u(z)\geq \bar{h}>0 \quad \text{on} \quad B^-_{\bar{\theta}r}.
\end{eqnarray}
\end{lemma}

{\bf Proof}: Let $h_0$ and $\theta$ be fixed such that the conclusions in Proposition 4.5 and Lemma 4.6 hold. Set
$$w(z)=ln^+ \left( \frac{h}{u+h^{\frac{9}{8}}}\right), \quad 0<h\leq h_0$$

It then follows from Proposition 4.5 that
\begin{eqnarray}
f_{B^-_{\theta r}} ((w-I_0)^+)^2dz\leq c\frac{\theta^2 r^2}{|B^-_{\theta r}|} \int_{B^-_{r/\theta}} |\partial_y w(z)|^2 dz.
\end{eqnarray}

Where $f_\Omega f\,dz$ denotes the average of $f$ on $\Omega$. To estimate the integral on the right hand side of (4.63), we note that $w(z)=ln^+(\frac{1}{\tilde{u}+h^{\frac{1}{8}}})$ with $\tilde{u}=\frac{u}{h}$ being a non-negative solution to (4.15). One then can follow the argument for (4.22) in the proof of Proposition 4.2 that for $\tau<t$ in $(-(\tilde{r})^2,0)$ with $\tilde{r}=\frac{r}{\theta}<1$, it holds that
\begin{eqnarray}
\begin{array}{cl}
& \displaystyle \int_{c^-_{\tilde{r}}} w(z)dxdy+\int^t_\tau\!\!\int_{c^-_{(1+\delta)\tilde{r}}}\chi^2_\delta a|\partial_y\,w|^2 dz\\
& \displaystyle \int_{c^-_{(1+\delta)\tilde{r}}} \chi^2_\delta w(l,y,\tau)dxdy - \int^t_\tau\!\!\int_{c^-_{(1+\delta)\tilde{r}}} 2\chi_\delta (\partial y\chi_\delta)a\partial_y w\,dz - \int^t_\tau\!\!\int_{c^-_{(1+\delta)\tilde{r}}} \chi^2_\delta\, y\, \partial_x\, w\, dz
\end{array}
\end{eqnarray}
where $\delta\in(0,1)$ so that $(1+\delta)\tilde{r}<1$, and $\chi_\delta=\chi_\delta(|y|)$ and $\chi_\delta\in c^\infty[0,\infty)$, $0<\chi_\delta(s)\leq 1$, $\chi_\delta(s)=1$ for $0\leq s\leq\tilde{r}$, $\chi_\delta(s)=0$ for $s>(1+\delta)\tilde{r}$, and
\begin{eqnarray}
|\chi'_\delta(s)|\leq \frac{2}{\delta\tilde{r}} \quad \text{for all} \quad s\geq 0
\end{eqnarray}

As in the derivation of (4.23), one has from (4.65) that
\begin{eqnarray*}
\begin{array}{rcl}
\displaystyle \int^t_\tau\!\!\int_{c^-_{(1+\delta)\tilde{r}}} 2\chi_\delta\partial_y\chi_\delta a\partial_y w dz & \leq & \displaystyle \frac{1}{2}\int^t_\tau\!\!\int_{c^-_{(1+\delta)\tilde{r}}} a\chi^2_\delta(\partial_y w|^2dz+\frac{16\Lambda}{\delta^2\tilde{r}^2}(t-\tau)|c^-_{(1+\delta)\tilde{r}}|\\
& \leq & \displaystyle \frac{1}{2}\int^t_\tau\!\!\int_{c^-_{(1+\delta)\tilde{r}}} a\chi^2_\delta|\partial_y w|^2dz+\frac{16\Lambda}{\delta^2}(1+\delta)^4\tilde{r}^4,\\
\displaystyle -\int^t_\tau\!\!\int_{c^-_{(1+\delta)\tilde{r}}} \chi^2_\delta y\partial_x w\,dz & = & \displaystyle -\int^t_\tau\!\!\int^{(1+\delta)\tilde{r}}_{-(1+\delta)\tilde{r}} \chi^2_\delta yw \bigg|^{((1+\delta)\tilde{r})^3}_{-((1+\delta)\tilde{r})^3} dyds\\
& \leq & \displaystyle (t-\tau)(1+\delta)^2 \tilde{r}^2 ln h^{-\frac{1}{8}}\leq (1+\delta)^2 \tilde{r}^4 ln\,h^{-\frac{1}{8}},
\end{array}
\end{eqnarray*}
and
\begin{eqnarray*}
\begin{array}{rcl}
\displaystyle \int_{c^-_{(1+\delta)\tilde{r}}} \chi^2_\delta w(x,y,\tau)dxdy & = & \displaystyle \int_{c^-_{(1+\delta)\tilde{r}}\cap\{\tilde{u}\leq 1\}} \chi^2_\delta w(x,y,\tau)dxdy\\
& \leq & \displaystyle ln\,h^{-\frac{1}{8}}\, \text{mes} \, c^-_{(1+\delta)\tilde{r}}=(1+\delta)^4 \tilde{r}^4 ln\,h^{-\frac{1}{8}}.
\end{array}
\end{eqnarray*}

These and (4.64) yield that for $\tau$ and $t$, $\tau<t$ in $(-\tilde{r},0)$,
$$\int^t_\tau\!\!\int_{c^-_{(1+\delta)\tilde{r}}} \chi^2_\delta a|\partial_y w(z)|^2dz\leq ((1+\delta)^4 (1+\frac{16\Lambda}{\delta^2}))\tilde{r}^4 ln\,h^{-\frac{1}{8}}.$$

This, together with the construction of $\chi_\delta$ and (4.16), implies that
\begin{eqnarray}
\int_{B^-_{\tilde{r}}} |\partial_y w(z)|^2dz\leq C(\Lambda,\delta)\tilde{r}^4 ln\,h^{-\frac{1}{8}}
\end{eqnarray}

It follows from (4.66) and (4.63) that there exists a positive constant $C=C(\Lambda,\theta,\delta)>0$ such that
\begin{eqnarray}
f_{B^-} ((w-I_0)^+)^2 dz\leq C\frac{\theta^2\,r^2}{|B^-_{\theta r}|} \int_{B^-_{\tilde{r}}} |\partial_y w|^2 dz\leq C\,ln\,h^{-\frac{1}{8}}.
\end{eqnarray}

Based on (4.67), one can modify the Moser's iteration method for $(w-I_0)^+$ as in [10] to show that there exists a uniform constant $k\in(0,1)$ such that
\begin{eqnarray}
\sup_{B^-_{k\theta r}} ((w-I_0)^+)^2 \leq C\,ln\,h^{-\frac{1}{8}}.
\end{eqnarray}

Set $\bar{\theta}=k\theta$. Then (4.68) and Lemma 4.6 imply that for all $z\in B^-_{\bar{\theta}r}$,
\begin{eqnarray}
w(z)\leq I_0+C(ln\,h^{-\frac{1}{8}})^2 \leq \lambda_0\,ln\,h^{-\frac{1}{8}} + C(ln\,h^{-\frac{1}{8}})^2.
\end{eqnarray}

Since $\displaystyle \lim_{h\rightarrow 0^+}\frac{\lambda_1\,ln\,h^{-\frac{1}{8}}+C(ln\,h^{-\frac{1}{8}})^{\frac{1}{2}}}{ln(\frac{1}{2h^{\frac{1}{8}}})}=\lambda_0<1$ due to Lemma 4.6, thus there exists a constant $h_2\leq \min(h_0,2^{-8})$ such that
$$\lambda_0\,ln\,h^{-\frac{1}{8}}_2+C(ln\,h_2^{-\frac{1}{8}})^{\frac{1}{2}}\leq ln(\frac{1}{2h^{\frac{1}{8}}_0}).$$

Consequently, it holds that
\begin{eqnarray}
\max_{B^-_{\bar{\theta}r}} ln^+ (\frac{h_2}{u+h_2^{\frac{9}{8}}})\leq ln(\frac{1}{2h^{\frac{1}{8}}_2}).
\end{eqnarray}

It follows from (4.70) and the choice of $h_2$ that
\begin{eqnarray}
\min_{B^-_{\bar{\theta}r}}\,u\geq h^{\frac{9}{8}}_2
\end{eqnarray}
which yields the desired estimate (4.62). Hence the proof of Lemma 4.7 is completed.

As an immediately consequence of Lemma 4.7, the following desired oscillation estimate holds.

\begin{proposition}
Let $u$ be a weak solution of (4.15) in $B^-_1$ and $\bar{\theta}$ and $r$ be given as in Lemma 4.7. Then there exists an uniform constant $\bar{\beta}$, $0<\bar{\beta}<1$, such that
\begin{eqnarray}
\text{Osc}_{B^-_{\bar{\theta}r}} \,u\leq \bar{\beta}\,\text{Osc}_{B^-_r}\,u
\end{eqnarray}
where $\text{Osc}_Q\,f$ denotes the oscillation of $f$ over $Q$ for any domain $\Omega$.
\end{proposition}

{\bf Proof}: Since $u$ is bounded, one can assume that $M\equiv\max_{B^-_r}\,u=-m\equiv-\min_{B^-_r}\,u$ without loss of generality since otherwise, one may consider $u-\frac{1}{2}(M+m)$. Thus both $1+\frac{u}{M}$ and $1-\frac{u}{M}$ are non-negative weak solutions to (4.15), and at least one of them satisfies the main assumption (4.17). We treat the case that $1-\frac{u}{M}$ satisfies (4.17). Then Lemma 4.7 can be applied to $(1-\frac{u}{M})$ to get

$$(1-\frac{u}{M})(z)\geq \bar{h} \quad \text{on}\quad B^-_{\bar{\theta}r},$$
where $\bar{h}\in(0,1)$ is a constant given in Lemma 4.7. Thus one gets that $\max_{B^-_{\bar{\theta}r}}\,u\leq (1-\bar{h})M$. Consequently,
\begin{eqnarray*}
\begin{array}{rcl}
\displaystyle \text{Osc}_{B^-_{\bar{\theta}r}}\,u & = & \displaystyle \max_{B^-_{\bar{\theta}r}}\,u-\min_{B^-_{\bar{\theta}r}}\,u\leq\max_{B^-_{\bar{\theta}r}}\, u+M\\
& \leq & \displaystyle (1-\bar{h})M+M=(1-\frac{\bar{h}}{2})2M=(1-\frac{\bar{h}}{2})\text{Osc}_{B^-_r}\,u,
\end{array}
\end{eqnarray*}
which yields the desired estimate (4.72) with $\bar{\beta}=(1-\frac{\bar{h}}{2})$. Hence the proof of Proposition (4.72) is completed. $\square$

It follows from the oscillation estimate (4.72) in Proposition 4.8 and the standard regularity arguments [16] that the following statement holds.

\begin{proposition} Let $u$ be a weak solution to (4.15) in $B^-_1$. Then there exists a positive constant $\delta>0$ such that $u$ is H\"{o}lder continuous on $B^-_\delta$ which is a small neighborhood of $z=0$.
\end{proposition}

Now, Proposition 4.9 shows that the weak solution $u$ to (4.2) is H\"{o}lder continuous in an interior point in $Q_T$. Then the standard regularity theory for ultra-parabolic equation [1, 14]. $u$ is in fact $C^\infty$ smooth in $Q_T$. Thus the proof of Theorem 4.1 is completed.

{\bf Remark:} In the preparation of the current paper, there are some new developments in H\"older continuity of weak solutions to a class of ultraparabolic equations, which generalized Theorem 4.1 [17,18,19,20,21]. In particular, in [19,21] their results include a different proof of Theorem 4.1 in the general case.

\end{document}